\documentclass[10pt,a4paper]{article}

\usepackage{amsmath,amsthm}
\usepackage[english]{babel}
\usepackage[latin1]{inputenc}
\usepackage{amsfonts}
\usepackage{amssymb}
\usepackage{amsmath}
\usepackage{graphics}
\usepackage{tabularx}

\usepackage[active]{srcltx}

\topmargin = - 2 true cm
\theoremstyle{plain}

\newtheorem{stat}{Statement}[section]
\newtheorem{thm}[stat]{Theorem}
\newtheorem{prop}[stat]{Proposition}
\newtheorem{coro}[stat]{Corollary}
\newtheorem{lemma}[stat]{Lemma}
\newtheorem{thm*}{Theorem C.0.\!\!}
\theoremstyle{definition}
\newtheorem{rem}[stat]{Remark}

\begin{document}


\newcommand{\N}{\mathbf{N}}
\newcommand{\Z}{\mathbf{Z}}
\newcommand{\Q}{\mathbb{Q}}
\newcommand{\R}{\mathbf{R}}
\newcommand{\C}{\mathbb{C}}

\renewcommand{\P}{{\rm P}}
\newcommand{\E}{{\rm E}}
\newcommand{\Filt}{\mathcal{F}}

\newcommand{\ind}{\boldsymbol{1}}
\newcommand{\Lip}{{\rm Lip}}
\newcommand{\Ell}{{\rm L}}
\newcommand{\const}{\text{\rm const}}

\newcommand{\espilon}{\epsilon}

\renewcommand{\geq}{\geqslant}
\renewcommand{\leq}{\leqslant}

\numberwithin{equation}{section}





\begin{center}
{\bf \Large Intermittency and chaos for a stochastic non-linear wave equation in dimension 1.}
\vskip 16pt

Daniel Conus\footnote[1]{Department of Mathematics, Lehigh University, daniel.conus@lehigh.edu.},
Mathew Joseph\footnote[2]{Department of Mathematics, University of Utah, joseph@math.utah.edu.},
Davar Khoshnevisan\footnote[3]{Department of Mathematics, University of Utah, davar@math.utah.edu.},
Shang-Yuan Shiu\footnote[4]{Academica Sinica, shiu@math.sinica.edu.tw.}.
\vskip 12pt 

November 2011

\end{center}

\begin{abstract}
We consider a non-linear stochastic wave equation driven by space-time white noise in dimension 1. First of all, we state some results about the intermittency of the solution, which had only been carefully studied in some particular cases so far. Then, we establish a comparison principle for the solution, following the ideas of Mueller. We think it is of particular interest to obtain such a result for an hyperbolic equation. Finally, using the results mentioned above, we aim to show that the solution exhibits a chaotic behavior, in a similar way as was established in \cite{CJK} for the heat equation. We study the two cases where 1. the initial conditions have compact support, where the global maximum of the solution remains bounded and 2. the initial conditions are bounded away from 0, where the global maximum is almost surely infinite. Interesting estimates are also provided on the behavior of the global maximum of the solution.
\end{abstract}
\vskip 2in
\vskip 12pt

%



\section{Introduction}\label{sec:intro}

We consider the family of Stochastic Partial Differential Equations given by
\begin{equation} \label{eqn:wave}
\frac{\partial^2 u}{\partial t^2}(t,x) = \kappa^2 \frac{\partial^2 u}{\partial x^2} u(t,x) + \sigma(u(t,x)) \dot{W}(t,x), \qquad t > 0, x \in \mathbf{R},
\end{equation}
where $\sigma:\mathbf{R} \rightarrow \mathbf{R}$ is a globally Lipschitz function with Lipschitz constant $\Lip_{\sigma}$, the noise $(\dot{W}(t,x), t \geq 0, x \in \mathbf{R})$ is space-time white noise and $\kappa > 0$. We consider non-random, bounded and measurable initial condition $u_0 : \mathbf{R} \rightarrow \mathbf{R}_+$ and initial derivative $v_0 : \mathbf{R} \rightarrow \mathbf{R}$.

Equation \eqref{eqn:wave} has been studied by Carmona and Nualart \cite{CarmonaNualart} and Walsh \cite{walsh}. There are also results available in the more delicate setting where
$x\in\R^d$ for $d>1$; see Conus and Dalang \cite{conus_dalang}, Dalang \cite{dalang}, Dalang and Frangos \cite{dalang_frangos}, and Dalang and Mueller \cite{dalang_mueller}.

The parabolic equivalent to equation \eqref{eqn:wave} is a well-studied family of Stochastic PDEs. In particular, it contains the well-known Parabolic Anderson Model. For more about this parabolic family, we refer the reader to Foondun and Khoshnevisan \cite{foondun_khoshnevisan}.

It is well-known that the Green function for the wave operator in spatial dimension 1 is 
\begin{equation} \label{eqn:green_wave}
	\Gamma_t(x) := \frac{1}{2} \ind_{[-\kappa t,\kappa t]}(x)\qquad\text{for $t>0$ and $x\in\R$}.
\end{equation}
According to the theories of Walsh \cite{walsh} and Dalang \cite{dalang}, the stochastic wave equation \eqref{eqn:wave} has an a.s.-unique mild solution $\{u(t,x), t > 0, x \in \R\}$, which satisfies
\begin{equation} \label{eqn:mild}
	u(t,x) = U_{0}(t,x) + V_{0}(t,x) + \int_{0}^{t} \int_{\R} \Gamma_{t-s}(y-x) \sigma(u(s,y))\, W(ds, dy),
\end{equation}
where the stochastic integral is defined in the sense of Walsh \cite{walsh} and
\begin{equation}
U_0(t,x) := \frac{1}{2} (u_0(x+\kappa t) + u_0(x-\kappa t)),
\end{equation}
and
\begin{equation}
V_0(t,x) := \frac{1}{2} \int_{x-\kappa t}^{x+\kappa t} v_0(y) \, dy.
\end{equation}

We remind that a process $\{u(t,x) : t > 0, x \in \mathbf{R}\}$ is \emph{weakly intermittent} if the upper Lyapunov exponents
\begin{equation} \label{eqn:lyapunov}
\bar{\gamma}(p) := \limsup_{t \rightarrow \infty} \frac{1}{t} \sup_{x \in \R} \log \E[|u(t,x)|^p] \qquad (1 \leq p < \infty)
\end{equation}
satisfy
$$\bar{\gamma}(2) > 0 \qquad \mbox{ and } \qquad \bar{\gamma}(p) < \infty, \quad \mbox{for all } p \geq 2.$$

A weakly intermittent process develops very high peaks as the time parameter gets large. Motivated by conjectures from physics related to the so-called KPZ equation, intermittency is of major importance in the parabolic equivalent family of \eqref{eqn:wave}, which includes the well-known Parabolic Anderson model. We refer to Foondun and Khoshnevisan \cite{foondun_khoshnevisan} for more information on the role of intermittency in the parabolic case.

As far as the wave equation is concerned, Dalang and Mueller \cite{dalang_mueller_2} have shown that the solution to \eqref{eqn:wave} is intermittent in the case that $u_0$ and $v_0$ are both constant functions, $\sigma(u) = \lambda u$ (Hyperbolic Anderson Model), and $\dot{W}$ is a colored noise, using a Feynman-Kac type representation. Intermittency in the case where $u_0$ and $v_0$ have compact support and $\dot{W}$ is a space-time white noise was established by Conus and Khoshnevisan in \cite{conus_khoshnevisan}, using a stochastic Young-type inequality. They also obtained results on the position of peaks of the solution.

The purpose of this paper is to carefully cover the study of intermittency and chaotic properties of the stochastic non-linear wave equation \eqref{eqn:wave}. We mainly follow the ideas of Conus, Joseph and Khoshnevisan \cite{CJK}, who establish precise estimates on the asymptotic behavior of $\sup_{|x|\leq R} u_t(x)$ as $R \rightarrow \infty$ for fixed $t > 0$ in the parabolic case. They thereby show that the solution to the stochastic heat equation exhibits a chaotic behavior. We would like to address a similar program in the case of the hyperbolic analogue family of SPDEs \eqref{eqn:wave}. We will obtain similar estimates on the supremum of the solution to the hyperbolic equation, showing that it also exhibits a chaotic behavior. The exact rates of the estimates are different than the ones of the parabolic case and the comparison is of particular interest, since it illustrates both similarities and discrepancies between the wave and the heat equation.

Since several procedures used below are similar to the ones used for the parabolic case and only differ in computations, we will assume that the reader is familiar with the ideas of Foondun and Khoshnevisan \cite{foondun_khoshnevisan, foondun_khoshnevisan_2} and Conus, Joseph and Khoshnevisan \cite{CJK} and we do not hesitate to refer to these papers whenever necessary. However, some results are specific to the hyperbolic case and we give complete proofs in those cases.

Section \ref{sec:inter} below is mainly a reminder about the intermittency results for \eqref{eqn:wave}. These results are mostly known but have never been stated in a uniform way, which we would like to do. We will derive a comparison principle for the stochastic wave equation in Section \ref{sec:comp}. This result will be needed in order to establish supremum estimates, although it is of specific interest in itself, since comparison principles for wave equations are pretty uncommon. Section \ref{sec:bdd} will be devoted to the study of the case where the initial condition has compact support. We will show that $\sup_{x\in\R} u(t,x) < \infty$ a.s. when the initial conditions have compact support. This result is intuitively a consequence of \cite{conus_khoshnevisan}, but needs to be formally established. Several moment and tail probability estimates for different behavior of the non-linearity $\sigma$, as well as localization properties are proved in Sections \ref{sec:moment_tail} and \ref{sec:loc}. Finally, these results will be used in Section \ref{sec:chaos}, where we will state and prove the main results of this paper about the asymptotic behavior of $\sup_{|x| \leq R} u_t(x)$.

In \cite{CJK}, a result on the behavior of the supremum in the case of the Anderson model is obtained, based on the sharp moment estimates of Bertini and Cancrini \cite{bertini_cancrini}. Such estimates are not available in the hyperbolic case. The hyperbolic Anderson model is subject of ongoing research. The study of intermittency and chaotic properties for the parabolic equation has also led to further developments, among which we can cite the study of the equation driven by colored-noise (see \cite{CJKS}) and the study of the size of the intermittent islands (see \cite{CJK2}). Similar extensions for the solution to the hyperbolic equation (which would relate to \cite{dalang_mueller_2}) is subject to ongoing research as well.

We would like to introduce some notation and preliminary results that are used throughout the paper. For a random variable $Z$, we denote by $\|Z\|_p := \E[|Z|^p]^{1/p}$ the standard norm on $L^p(\Omega)$ ($p \geq 1$). On several occasions, we apply the Burkholder--Davis--Gundy inequality \cite{Burkholder,BDG,BG} for continuous $L^2(\Omega)$ martingales:
If $\{X_t\}_{t \geq 0}$ is a continuous $L^2(\Omega)$ martingale with running maximum $X^*_t:=\sup_{s\in[0,t]}|X_s|$ and quadratic variation process $\langle X\rangle$, then for all 
real numbers $p \geq 2$ and $t>0$,
\begin{equation} \label{eqn:BDG}
\left\| X_t^*\right\|_p \leq \left\| 4p\langle X\rangle_t \right\|_{p/2}^{1/2}.
\end{equation}
The factor $4p$ is the asymptotically-optimal bound of Carlen and Kree \cite{CK} for the sharp constant in the Burkholder--Davis--Gundy
inequality that is due to Davis \cite{Davis}. We will also sometimes use the notation
\begin{equation}\label{eqn:u_bar}
\underline{u}_0:=\inf_{x\in \R^d}u_0(x),\quad \overline{u}_0:=\sup_{x\in \R^d} u_0(x), \quad \underline{v}_0:=\inf_{x\in \R^d}v_0(x),\qquad \overline{v}_0:=\sup_{x\in \R^d} v_0(x).
\end{equation}
Further, for a random field $\{Z(t,x) : t > 0, x \in \mathbf{R}\}$, we denote by $\|\cdot\|_{p,\beta}$ the norm defined by
\begin{equation}
\|Z\|_{p,\beta} := \left(\sup_{t \geq 0} \sup_{x \in \R} e^{-\beta t} \E[|Z(t,x)|^p] \right)^{1/p}.
\end{equation}
This norm is the one used in \cite{foondun_khoshnevisan}. In the literature (see e.g. \cite{CJK}), we sometimes find similar results expressed with the norm
\begin{equation}
\mathcal{N}_{p,\beta}(Z) := \left(\sup_{t \geq 0} \sup_{x \in \R} e^{-\beta t} \|Z\|^2_p \right)^{1/2}.
\end{equation}
Although they are formally different, the two norms are equivalent and are related by
\begin{equation}
\mathcal{N}_{p,\beta}(Z) = \|Z\|_{p,\frac{p\beta}{2}}, \quad \|Z\|_{p,\beta} = \mathcal{N}_{p,\frac{2\beta}{p}}(Z)
\end{equation}
Even though this detail does not matter when it comes to proving that one norm or the other is finite, it is relevant when we want to keep track of the correct constants. Finally, we have the following immediate results about the Green function $\Gamma$:
\begin{equation} \label{eqn:gamma_L2}
\|\Gamma_t\|^2_{L^2(\R)} = \frac{\kappa t}{2},
\end{equation}
\begin{equation} \label{eqn:gamma_L2b}
\int_{0}^{t} \|\Gamma_s\|^2_{L^2(\R)} \, ds = \frac{\kappa t^2}{4},
\end{equation}
and
\begin{equation} \label{eqn:upsilon}
\Upsilon(\beta) := \int_{0}^{\infty} e^{-\beta t} \|\Gamma_t\|^2_{L^2(\R)} \, dt = \frac{\kappa}{2 \beta^2}.
\end{equation}
The last expression is labelled $\Upsilon(\beta)$ by analogy with the potential theory for L\'evy processes (see \cite{foondun_khoshnevisan}).


\section{Intermittency} \label{sec:inter}

We are now ready to state and prove intermittency for the solution to \eqref{eqn:wave} under general assumptions on the initial condition and the non-linearity $\sigma$ in the case where $\dot{W}$ is space-time white noise. These results follow closely the techniques of Foondun and Khoshnevisan \cite{foondun_khoshnevisan} for the hest equation and of Conus and Khohsnevisan \cite{conus_khoshnevisan} for the wave equation with compact support initial data. In particular, we will use the stochastic Young-type inequality of \cite{conus_khoshnevisan}.

In order to establish weak intermittency of the solution to \eqref{eqn:wave}, we need to obtain two different results. We need to prove an upper bound in that $\bar{\gamma}(p) < \infty$ for all $p \geq 1$ and, then, we need to prove that $\bar{\gamma}(2) > 0$, where $\bar{\gamma}(p)$ is defined in \eqref{eqn:lyapunov}.

\begin{thm} \label{thm:inter_up}
Let $u$ denote the solution to \eqref{eqn:wave} and $\bar{\gamma}(p)$ be defined by \eqref{eqn:lyapunov}. Then, for all $p \geq 2$,
$$\bar{\gamma}(p) \leq p^{3/2} \Lip_{\sigma} \sqrt{\kappa/2}.$$
\end{thm}

\begin{rem} \label{rem:inter_up}
A careful look at the proof of Theorem \ref{thm:inter_up} shows that the factor $2^{3/2}$ is not needed in the case where $p=2$. Indeed, the optimal constant in the Burkholder-Davis-Gundy inequality is 1 when $p=2$.
\end{rem}

For the next result, we let 
\begin{equation} \label{eqn:ell}
\Ell_{\sigma} := \inf_{x \neq 0} \left|\frac{\sigma(x)}{x}\right|.
\end{equation}

\begin{thm} \label{thm:inter_low}
If $\underline{u}_0 = \inf_{x \in \R} u_0(x) > 0$, $v_0 \geq 0$ and $\Ell_{\sigma} > 0$, then
$$\bar{\gamma}(2) \geq \Ell_{\sigma} \sqrt{\kappa/2}$$
\end{thm}

Theorems \ref{thm:inter_up} and \ref{thm:inter_low} are similar to Theorems 2.1 and 2.7 of \cite{foondun_khoshnevisan}. Together, they prove that the solution $u$ is weakly intermittent provided that $u_0$ is bounded away from 0, $v_0 \geq 0$ and $\sigma$ has linear growth. Intermittency in the case where $u_0$ and $v_0$ have compact support has been studied in \cite{conus_khoshnevisan}, see also Section \ref{sec:bdd}. Linear growth of $\sigma$ is somehow necessary for intermittency as the following result suggests (see \cite[Thm.2.3]{foondun_khoshnevisan} for the parabolic case).

\begin{thm} \label{thm:no_inter}
If $u_0, v_0$ and $\sigma$ are all bounded functions, then for all $p \geq 2$,
$$\E[|u(t,x)|^p] = o(t^p), \qquad \mbox{ as } t \rightarrow \infty.$$
\end{thm}

These result are known, mostly from \cite{foondun_khoshnevisan} and \cite{conus_khoshnevisan}, but for the sake of completeness and since they were never really stated independently, we give an outline of their proof. We will use a stochastic Young-type inequality for stochastic convolutions (Proposition \ref{prop:young} below), which is a direct consequence of Proposition 2.5 of \cite{conus_khoshnevisan} (we omit the proof).

For a random-field $(Z(t,x),t>0,x\in\R)$, we denote by $\Gamma*Z\dot{W}$ the random-field defined by
$$(\Gamma*Z\dot{W})(t,x) = \int_{0}^{t} \int_{\R} \Gamma_{t-s}(y-x) Z(s,y) W(ds,dy),$$
when the stochastic integral is well-defined in the sense of Walsh \cite{walsh}.

\begin{prop} \label{prop:young}
For all $\beta > 0$ and $p \geq 2$,
\begin{equation}
\|\Gamma*Z\dot{W}\|_{p,\beta} \leq 2 \sqrt{p} \, \Upsilon\left(\frac{2\beta}{p}\right)^{1/2} \|Z\|_{p,\beta}.
\end{equation}
\end{prop}

We are now ready to prove the main results of this section.\\

\noindent{\bf Proof of Theorem \ref{thm:inter_up}:}
Since $u_0$ and $v_0$ are bounded above by $\overline{u}_0$ and $\overline{v}_0$ respectively, we clearly have
\begin{equation} \label{eqn:initial}
\sup_{x\in\R} |U_0(t,x) + V_0(t,x)| \leq \overline{u}_0 + \overline{v}_0 \kappa t
\end{equation}
Now, using \eqref{eqn:mild}, \eqref{eqn:initial}, Proposition \ref{prop:young}, and the fact that $|\sigma(u)| \leq |\sigma(0)| + \Lip_{\sigma}|u|$, we have for all $\beta > 0$ and $p \geq 2$,
\begin{equation}
\|u\|_{p,\beta} \leq \overline{u}_0 + \frac{p^2\kappa}{\beta^2} \,\overline{v}_0 + 2\sqrt{p} \, \Upsilon\left(\frac{2\beta}{p}\right)^{1/2} (|\sigma(0)| + \Lip_{\sigma} \|u\|_{p,\beta}).
\end{equation}
This shows that $\|u\|_{p,\beta} < \infty$ provided that $2\sqrt{p} \, \Lip_{\sigma} \Upsilon(2\beta/p)^{1/2} < 1$, which is equivalent to
\begin{equation}
\beta > p^{3/2} \Lip_{\sigma} \sqrt{\kappa/2}.
\end{equation}
Theorem \ref{thm:inter_up} follows.
\hfill $\blacksquare$ \\

\noindent{\bf Proof of Theorem \ref{thm:inter_low}:}
This proof follows precisely the ideas of the proof of Theorem 2.7 of \cite{foondun_khoshnevisan} which directly apply here. If we first prove that
\begin{equation}\label{eqn:cond_infinite}
\int_{0}^{\infty} e^{-\beta t} \E[|u(t,x)|^2] \, dt = \infty \quad \mbox{ provided that } \Upsilon(\beta) > \Ell_{\sigma}^{-2}, 
\end{equation}
then Theorem \ref{thm:inter_low} follows. We prove \eqref{eqn:cond_infinite} exactly as in the proof of Theorem 2.7 of \cite{foondun_khoshnevisan} using Laplace transforms and a renewal equation argument. The only difference is that we need to check that
\begin{equation}
\int_{0}^{\infty} e^{-\beta t} (U_0(t,x) + V_0(t,x)) \, dt > 0,
\end{equation}
which is a direct consequence of the fact that $V_0(t,x) \geq 0$ and $U_0(t,x) > \underline{u}_0 > 0$ for all $t > 0$, under the assumption of Theorem \ref{thm:inter_low}.
\hfill $\blacksquare$ \\

\noindent{\bf Proof of Theorem \ref{thm:no_inter}:}
Since $\sigma$ is bounded, a direct consequence of \eqref{eqn:mild}, \eqref{eqn:initial}, the BDG inequality, and \eqref{eqn:gamma_L2b} is that
\begin{eqnarray}
\|u(t,x)\|_p & \leq & \overline{u}_0 + \overline{v}_0 \kappa t + 2 \sqrt{p} \left(\sup_{x \in \R} \sigma(x)\right) \left(\int_{0}^{t} \|p_s\|^2_{L^2{\R}} \, ds\right)^{1/2} \nonumber \\
& \leq & \overline{u}_0 + \overline{v}_0 \kappa t + \sqrt{p\kappa} \left(\sup_{x \in \R} \sigma(x)\right) t. \label{eqn:beta_large}
\end{eqnarray}
Theorem \ref{thm:no_inter} follows.
\hfill $\blacksquare$


\section{Comparison principle} \label{sec:comp}

The purpose of this section is to establish a comparison principle for the stochastic wave equation. The main idea behind the proof is to approximate the solution $u$ of \eqref{eqn:wave} by the solution to an infinite system of Stochastic Differential Equations, which would themselves satisfy a comparison principle. The ideas of this result are similar to the ones developped by Mueller \cite{mueller} for the parabolic equation. Without loss of generality, for this section, we will only consider the case $\kappa=1$.

\begin{thm} \label{thm:comp_p}
Let $u^{(i)}$ denote the solution to \eqref{eqn:wave} with respectively $u^{(i)}_0$ and $v^{(i)}_0$ as initial condition and initial derivative ($i=1,2$). If 
\begin{equation} \label{eqn:comp_p_assumption}
u^{(1)}_0(x) \geq u^{(2)}_0(x) \qquad \mbox{ and } \qquad v^{(1)}_0(x) \geq v^{(2)}_0(x), \qquad \mbox{ for all } x \in \R,
\end{equation}
and $u^{(1)}_0$ and $u^{(2)}_0$ are Hölder continuous of order $\alpha \in (0,1]$, then for all $t >0$ and $x \in \R$,
\begin{equation}
u^{(1)}(t,x) \geq u^{(2)}(t,x) \qquad \mbox{ a.s.}
\end{equation}
\end{thm}

Before we turn to the proof of Theorem \ref{thm:comp_p} itself, we prove a series of results about a discrete space-time approximation of the solution $u$.

First of all, we define a dyadic approximation of space and time parameters. Let $n \in \mathbf{N}$. For every $x \in \R$, we write
\begin{equation} \label{eqn:approx_x}
\phi_n(x) := \frac{k}{2^n} \qquad \mbox{ if } \qquad \frac{k}{2^n} - \frac{1}{2^{n+1}} \leq x < \frac{k}{2^n} + \frac{1}{2^{n+1}} \quad (k \in \mathbf{Z}),
\end{equation}
and for every $t \geq 0$, we write
\begin{equation} \label{eqn:approx_t}
\psi_n(t) := \frac{m}{2^n} \qquad \mbox{ if } \qquad \frac{m}{2^n} \leq t < \frac{m+1}{2^n}, \quad (m \in \mathbf{N}).
\end{equation}
Clearly, 
\begin{equation} \label{eqn:xt_approx}
|x-\phi_n(x)| < 2^{-n} \qquad \mbox{ and } \qquad |t - \psi_n(t)| < 2^{-n}.
\end{equation}

Now, we would like to build a dyadic approximation to the Green function $\Gamma$. For $k \in \mathbf{Z}$, $m \in \mathbf{N}$, we set
\begin{equation} \label{eqn:gamma_n}
\Gamma^{(n)}_{m/2^n}\left(\frac{k}{2^n}\right) := \left\{
\begin{array}{ll}
1 & \mbox{ if } k < m, \\
1/2 & \mbox{ if } k=m, \\
0 & \mbox{ if } k > m.
\end{array}\right.
\end{equation}
The approximation \ref{eqn:gamma_n} actually corresponds to an $L^1$-type approximation, where we approximate the function by averaging over dyadic intervals. This explains the $1/2$ for the case $k=m$.

Since we would like to approximate the integrand of the stochastic integral, namely $\Gamma_{t-s}(y-x)$, we define a function $\Gamma_n$ by
\begin{equation}
\Gamma_n(t,s;x,y) := \Gamma^{(n)}_{\psi_n(t)-\psi_n(s)}(\phi_n(y)-\phi_n(x)),
\end{equation}
for $0 \leq s \leq t$ and $x,y\in\R$, where $\Gamma^{(n)}$ is the function defined in \eqref{eqn:gamma_n}.

In comparing with standard approximations for the parabolic Green function, notice that we do not need the time mesh to be of the order of the square root of the space mesh. The space-time behavior of the hyperbolic Green function $\Gamma$ allows to consider the same mesh for both. We choose the approximation \eqref{eqn:approx_x} for the space variable in order to preserve symmetry.

Finally, we are able to build our approximation of the solution $u$. For $n \in \mathbf{N}$, we define the random field $(u_n(t,x), t > 0, x \in \R)$ to be the solution to
\begin{equation} \label{eqn:approx}
u_n(t,x) = U_0(\psi_n(t),\phi_n(x)) + V_0(\psi_n(t), \phi_n(x)) + \int_{0}^{t} \int_{\R} \Gamma_n(t,s;x,y) \sigma(u_n(s,y)) W(ds,dy).
\end{equation}

\begin{prop} \label{prop:approx_ex}
The sequence of random fields $\{u_n(t,x), t>0,x\in\R\}_{n\in\N}$ is well-defined and satisfies
\begin{equation}\label{eqn:mom_unif}
\sup_{n\in\N}\sup_{0\leq t \leq T} \sup_{x \in \R} \E[|u_n(t,x)|^2] < \infty.
\end{equation}
\end{prop}

\noindent{\bf Proof.}
Existence and uniqueness of each of the $u_n$ is a direct consequence of the existence and uniqueness result for the original equation \eqref{eqn:wave}, provided
\begin{equation} \label{eqn:Dalang_n}
\sup_{t \in [0,T]} \sup_{x \in \R} \int_{0}^{t} \|\Gamma_n(t,s;x,\cdot)\|^2_{L^2(\R)} \, ds < \infty.
\end{equation}
In addition, we directly obtain \eqref{eqn:mom_unif} if we prove that \eqref{eqn:Dalang_n} holds uniformly in $n$.

For $t \geq 0$ such that $\psi_n(t) = m 2^{-n}$, we have
\begin{eqnarray}
\int_{0}^{t} \left\|\Gamma_n\left(t,s;x,\cdot\right)\right\|^2_{L^2(R)} \, ds & = &  \int_{0}^{t} ds \int_{\R} dy \, \Gamma^{(n)}_{\psi_n(t)-\psi_n(s)}(\phi_n(y)-\phi_n(x)) \\
& = &  \int_{0}^{t} ds \int_{\R} dy \, \Gamma^{(n)}_{\psi_n(t)-\psi_n(s)}(\phi_n(y)) \\
& \leq &  \int_{0}^{t} \left(2(\psi_n(t)-\psi_n(s)) + 2^{-n}\right) \, ds \\
& = & \int_{0}^{\psi_n(t)} 2(\psi_n(t)-\psi_n(s))\, ds + \frac{t}{2^n} \\
& = & \sum_{k=0}^{m-1} 2^{-n} \left(\frac{2(m-k)}{2^n}\right) + \frac{t}{2^n} \\
& = & \frac{m(m-1)}{4^n} + \frac{t}{2^n} \leq t^2 + t.
\end{eqnarray}
This proves \eqref{eqn:Dalang_n} and Proposition \ref{prop:approx_ex}.
\hfill $\blacksquare$ \\

Now, we would like to prove that $u_n$ constitutes an approximation to $u$.

\begin{prop} \label{prop:conv}
Let $u$ be the solution to \eqref{eqn:wave} with $u_0$ Hölder-continuous of order $\alpha \in (0,1]$ and let $u_n$ be defined by \eqref{eqn:approx}. Then, for every $T > 0$,
\begin{equation}
\sup_{0 \leq t\leq T} \sup_{x \in \mathbf{R}} \|u_n(t,x)-u(t,x)\|_2 \longrightarrow 0, \quad \mbox{ as } n \rightarrow \infty.
\end{equation}
\end{prop}

We will need a series of Lemmas.

\begin{lemma} \label{lem:estimate_1}
Let $t = \frac{m}{2^n}$. Then,
\begin{equation}
\int_{0}^{t} ds \int_{\R} dy \, |\Gamma^{(n)}_{\psi_n(t)-\psi_n(s)}(\phi_n(y))-\Gamma_{t-s}(y)|^2 \leq \text{\rm const} \cdot \frac{t}{2^n}.
\end{equation}
\end{lemma}

\noindent{\bf Proof.}
Since $t = \psi_n(t) = m/2^n$, and since $t-\psi_n(s) = \psi_n(t-s)+2^{-n}$ for every $0 \leq s \leq t$, a change of variable gives
\begin{eqnarray}
\lefteqn{\int_{0}^{t} ds \int_{\R} dy \, |\Gamma^{(n)}_{\psi_n(t)-\psi_n(s)}(\phi_n(y))-\Gamma_{t-s}(y)|^2} \nonumber \\
& = & \int_{0}^{m/2^n} ds \int_{\R} dy \, |\Gamma^{(n)}_{\psi_n(s)+2^{-n}}(\phi_n(y))-\Gamma_{s}(y)|^2 \nonumber \\
& = & \sum_{j=0}^{m-1} \sum_{k=-\infty}^{\infty} \int_{j/2^n}^{(j+1)/2^n} ds \int_{k/2^n-1/2^{n+1}}^{k/2^n+1/2^{n+1}} dy \, \left|\Gamma^{(n)}_{(j+1)/2^n}\left(\frac{k}{2^n}\right)-\Gamma_s(y)\right|^2 \nonumber \\
& = & 2 \sum_{j=0}^{m-1} \sum_{k=1}^{j} \int_{j/2^n}^{(j+1)/2^n} ds \int_{k/2^n-1/2^{n+1}}^{k/2^n+1/2^{n+1}} dy \, |1-\ind_{[-s,s]}(y)|^2 \nonumber \\
&& \qquad + 2 \sum_{j=0}^{m-1} \sum_{k=j+2}^{\infty} \int_{j/2^n}^{(j+1)/2^n} ds \int_{k/2^n-1/2^{n+1}}^{k/2^n+1/2^{n+1}} dy \, |\ind_{[-s,s]}(y)|^2 \nonumber \\
&& \qquad + 2 \sum_{j=0}^{m-1} \int_{j/2^n}^{(j+1)/2^n} ds \int_{j/2^n+1/2^{n+1}}^{j/2^n+3/2^{n+1}} dy \, \left|\frac{1}{2}-\ind_{[-s,s]}(y)\right|^2. \label{eqn:estimate_1}
\end{eqnarray}
Now, the first term does not vanish only if $k=j$, in which case it is bounded by $m/4^n$. The second sum is identically zero and the third sum is also bounded above by $m/(2\cdot4^n)$. Hence, the left-hand side of \eqref{eqn:estimate_1} is bounded above by $2m/4^n \leq 2t/2^n$. The result is proved.
\hfill $\blacksquare$ \\

\begin{lemma} \label{lem:estimate_3}
Let $t > 0$, then
\begin{equation}
\int_{0}^{t} ds \int_{\R} dy \, |\Gamma^{(n)}_{\psi_n(t)-\psi_n(s)}(\phi_n(y))-\Gamma_{t-s}(y)|^2 \leq \text{\rm const} \cdot \frac{t}{2^n}.
\end{equation}
\end{lemma}

\noindent{\bf Proof.}
If $t=\frac{m}{2^n}$, this was proved in Lemma \ref{lem:estimate_1}. Now, let $t \neq \psi_n(t) = \frac{m}{2^n}$. We have
\begin{eqnarray}
\lefteqn{\int_{0}^{m/2^n} ds \int_{\R} dy \, |\Gamma_{\psi_n(t)-s}(y)-\Gamma_{t-s}(y)|^2} \nonumber \\
& = & \int_{0}^{m/2^n} ds \int_{\R} dy \, |\Gamma_{t-\psi_n(t)+s}(y)-\Gamma_{s}(y)|^2 \nonumber \\
& \leq & \frac{m}{2^n} \left|t - \frac{m}{2^n}\right| \leq \frac{t}{2^n}. \label{eqn:estimate_3a}
\end{eqnarray}
Moreover, since $\psi_n(s) = \psi_n(t)$ whenever $\psi_n(t) \leq s \leq t$,
\begin{eqnarray}
\lefteqn{\int_{m/2^n}^{t} ds \int_{\R} dy \, |\Gamma^{(n)}_{\psi_n(t)-\psi_n(s)}(\phi_n(y))-\Gamma_{t-s}(y)|^2} \nonumber \\
& = & \int_{m/2^n}^{t} ds \int_{\R} dy \, \left|\frac{1}{2}\ind_{[-\frac{1}{2^{n+1}},\frac{1}{2^{n+1}}]}(\phi_n(y))-\ind_{[-(t-s),t-s]}(y)\right|^2 \nonumber \\
& = & \int_{m/2^n}^{t} ds \int_{-1/2^{n+1}}^{1/2^{n+1}} dy \, \left|\frac{1}{2}-\ind_{[-(t-s),t-s]}(y)\right|^2 \nonumber \\
& \leq & \frac{1}{4} \cdot \frac{1}{2^n} \left|t-\frac{m}{2^n}\right| \leq \frac{1}{4 \cdot 4^n}. \label{eqn:estimate_3b}
\end{eqnarray}
The result is a consequence of \eqref{eqn:estimate_3a}, \eqref{eqn:estimate_3b} and Lemma \ref{lem:estimate_1}.
\hfill $\blacksquare$ \\

\begin{lemma} \label{lem:estimate_2}
Let $t>0$ and $x \in \R$. Then,
\begin{equation}
\int_{0}^{t} ds \int_{\R} dy \, |\Gamma^{(n)}_{\psi_n(t)-\psi_n(s)}(\phi_n(y)-\phi_n(x))-\Gamma_{t-s}(y-x)|^2 \leq \text{\rm const} \cdot \frac{t}{2^n}.
\end{equation}
\end{lemma} 

\noindent{\bf Proof.}
If $x = \frac{k}{2^n}$, the result is an immediate consequence of Lemma \ref{lem:estimate_3} after a change of variable. Now assume that $x \neq \phi_n(x) = \frac{k}{2^n}$. Then,
\begin{eqnarray}
\lefteqn{\int_{0}^{t} ds \int_{\R} dy \, \left|\Gamma_{t-s}\left(y-\frac{k}{2^n}\right)-\Gamma_{t-s}(y-x)\right|^2} \nonumber \\
& = & \int_{0}^{t} ds \int_{\R} dy \, \left|\Gamma_{s}\left(y+x-\frac{k}{2^n}\right)-\Gamma_{s}(y)\right|^2 \nonumber \\
& = & \int_{0}^{t} ds \int_{\R} dy \, \left|\ind_{[-s,s]}\left(y+x-\frac{k}{2^n}\right)-\ind_{[-s,s]}(y)\right|^2 \leq t \cdot \left|x-\frac{k}{2^n}\right| \leq \frac{t}{2^n}. \label{eqn:estimate_2}
\end{eqnarray}
The result is a consequence of \eqref{eqn:estimate_2} and Lemma \ref{lem:estimate_2}.
\hfill $\blacksquare$ \\

We are now ready to prove Proposition \ref{prop:conv}. \\

\noindent{\bf Proof of Proposition \ref{prop:conv}:}
Since $u_0$ is assumed to be Hölder continuous of order $\alpha \in (0,1]$ and by \eqref{eqn:xt_approx}, we clearly have
\begin{equation}
\sup_{0\leq t \leq T} \sup_{x \in \R} |U_0(\psi_n(t),\phi_n(x)) - U_0(t,x)| \leq \const \cdot 2^{-n\alpha}.
\end{equation}
Moreover, since $v_0$ is bounded,
\begin{equation}
\sup_{0\leq t \leq T} \sup_{x \in \R} |V_0(\psi_n(t),\phi_n(x)) - V_0(t,x)| \leq \const \cdot 2^{-n}.
\end{equation}
Now, by \eqref{eqn:approx} and since stochastic integrals have zero expectation,
\begin{eqnarray}
\|u_n(t,x)-u(t,x)\|_2 & \leq & |U_0(\psi_n(t),\phi_n(x)) - U_0(t,x)| + |V_0(\psi_n(t),\phi_n(x)) - V_0(t,x)| \nonumber \\
&& \qquad + \; \|(\Gamma_n*\sigma(u_n)\dot{W})(t,x)-(\Gamma*\sigma(u)\dot{W})(t,x)\|_2 \nonumber \\
& \leq & \const \cdot 2^{-n\alpha} +  \|((\Gamma_n-\Gamma)*\sigma(u_n)\dot{W})(t,x)\|_2 \nonumber \\
&& \qquad + \; \|(\Gamma*(\sigma(u_n)-\sigma(u))\dot{W})(t,x)\|_2 \label{eqn:conv_a}
\end{eqnarray}
Standard estimates on the second moment of stochastic convolutions show that,
\begin{eqnarray}
\lefteqn{\|((\Gamma_n-\Gamma)*\sigma(u_n)\dot{W})(t,x)\|^2_2} \nonumber \\ 
& \leq & \left(\sup_{n \in \mathbf{N}} \sup_{0 \leq s \leq t} \sup_{x \in \R} \|\sigma(u_n(t,x))\|^2_2\right) \int_{0}^{t} ds \int_{\R} dy \, |\Gamma_n(t,s;x,y)-\Gamma_{t-s}(y-x)|^2 \nonumber \\
& \leq & \const \cdot \frac{t}{2^n}, \label{eqn:conv_b}
\end{eqnarray}
by Lemma \ref{lem:estimate_2} and Proposition \ref{prop:approx_ex}. Moreover,
\begin{eqnarray}
\|(\Gamma*(\sigma(u_n)-\sigma(u))\dot{W})(t,x)\|^2_2  & \leq & \int_{0}^{t} ds \int_{\R} dy \, \Gamma^2_{t-s}(y-x) \E[|\sigma(u_n(s,y)-\sigma(u(s,y)|^2] \nonumber \\
& \leq & \Lip^2_{\sigma} \int_{0}^{t} ds \, \left(\sup_{y \in \R} \E[|u_n(s,y)-u(s,y)|^2]\right) \, \|\Gamma_{t-s}\|^2_{L^2(\R)} \nonumber \\
& \leq & \Lip^2_{\sigma} \frac{t}{2} \int_{0}^{t} ds \, \sup_{y \in \R} \E[|u_n(s,y)-u(s,y)|^2]. \label{eqn:conv_c}
\end{eqnarray}
Now, setting $H_n(t):= \sup_{0 \leq s \leq t} \sup_{x \in \R} \E[|u_n(s,y)-u(s,y)|^2]$, \eqref{eqn:conv_a}, \eqref{eqn:conv_b} and \eqref{eqn:conv_c} show that
\begin{equation} \label{eqn:conv_gronwall}
H_n(t) \leq \const \cdot \left( 2^{-2n\alpha} + \frac{t}{2^n} + t \cdot \int_{0}^{t} H_n(s) \, ds \right).
\end{equation}
Gronwall's inequality applied to \eqref{eqn:conv_gronwall} implies that
\begin{equation}
\sup_{0 \leq t \leq T} H_n(t) \leq \const \cdot 2^{-2n(\alpha\wedge\frac{1}{2})},
\end{equation}
which proves the result.
\hfill $\blacksquare$ \\

We are now ready to prove the main result of this section. \\

\noindent{\bf Proof of Theorem \ref{thm:comp_p}:}
For $i=1,2$, let
\begin{equation}
U^{(i)}_0(t,x) := \frac{1}{2} \left(u^{(i)}_0(x+\kappa t) + u^{(i)}_0(x-\kappa t)\right),
\end{equation}
and
\begin{equation}
V^{(i)}_0(t,x) := \frac{1}{2} \int_{x-\kappa t}^{x+\kappa t} v^{(i)}_0(y) \, dy,
\end{equation}
and let $u^{(i)}_n(t,x)$ be defined according to \eqref{eqn:approx}, respectively with $U^{(i)}$ and $V^{(i)}$. First of all, we notice that \eqref{eqn:comp_p_assumption} implies that
\begin{equation} \label{eqn:comp_p_assumption_b}
U^{(1)}_0(t,x) \geq U^{(2)}_0(t,x) \qquad \mbox{ and } \qquad V^{(1)}_0(t,x) \geq V^{(2)}_0(t,x), \qquad \mbox{ for all } t > 0, x \in \R.
\end{equation}
Also notice that by construction, $u^{(i)}_n(t,x) = u^{(i)}_n(t,\phi_n(x))$ and, hence, is constant by intervals. This and \eqref{eqn:comp_p_assumption_b} together show that
\begin{equation} \label{eqn:comp_initial}
u^{(1)}_n(0,x) \geq u^{(2)}_n(0,x)
\end{equation}
for all $x \in \R$. Then, let $t > 0$, such that $\psi_n(t) = m/2^n$. Then, for each $i=1,2$ and $k \in \mathbf{Z}$,
\begin{eqnarray}
u^{(i)}_n\left(t,\frac{k}{2^n}\right) - u^{(i)}_n\left(\frac{m}{2^n},\frac{k}{2^n}\right) & = & \int_{m/2^n}^{t} \int_{\R} \Gamma^{(n)}_{0}\left(\phi_n(y)-\frac{k}{2^n}\right) \sigma(u^{(i)}_n(s,y)) W(ds,dy) \nonumber \\
& = & \frac{1}{2} \int_{m/2^n}^{t} \int_{k/2^n - 1/2^{n+1}}^{k/2^n + 1/2^{n+1}} \sigma(u^{(i)}_n(s,y)) W(ds,dy) \nonumber \\
& = & \frac{1}{2} \int_{m/2^n}^{t} \int_{k/2^n - 1/2^{n+1}}^{k/2^n + 1/2^{n+1}} \sigma\left(u^{(i)}_n\left(s,\frac{k}{2^n}\right)\right) W(ds,dy) \nonumber \\
\end{eqnarray}
Hence, if we set, for $i=1,2$, $n \in \mathbf{N}$ and $k \in \mathbf{Z}$,  
\begin{equation}
X^{(i)}_{n,k}(t) := u^{(i)}_n(t,k/2^n) \qquad \mbox{ and } \qquad B^{(n)}_k(t) = 2^n \, \int_{0}^{t} \int_{k/2^n - 1/2^{n+1}}^{k/2^n + 1/2^{n+1}} W(ds,dy),
\end{equation}
then $B^{(n)}_k(\cdot)$ is a Brownian motion and, for $t > 0$ such that $\psi_n(t) = m/2^n$,
\begin{equation} \label{eqn:SDE}
dX^{(i)}_{n,k}(t) := 2^{-(n+1)} \sigma(X^{(i)}_{n,k}(t)) \, dB^{(n)}_k(t)
\end{equation}
Hence, \eqref{eqn:SDE}, the comparison principle for Stochastic Differential Equations (see \cite[?]{revuz-yor}) and an induction argument on $m$ show that $u^{(1)}_n(t,x) \geq u^{(2)}_n(t,x)$ almost surely for all $t > 0$ and $x \in \R$. This and Proposition \ref{prop:conv} together prove the result.
\hfill $\blacksquare$ \\

\begin{rem}
Essentially, Theorem \ref{thm:comp_p} shows that a comparison principle holds for the stochastic wave equation because it can be seen as an infinite-dimensional system of SDEs, similarly as the parabolic equation. This holds even though the hyperbolic and parabolic are strongly different when seen from a PDE point-of-view.
\end{rem}


\section{Compact-support initial data} \label{sec:bdd}

This section is devoted to the study of the behavior of the supremum $\sup_{x\in\R} u(t,x)$ of the solution $u$ to \eqref{eqn:wave} when $t$ is fixed, in the case where the initial conditions $u_0$ and $v_0$ have compact support. We follow mainly the ideas of Foondun and Khoshnevisan \cite{foondun_khoshnevisan_2}. The matter is significantly easier in the hyperbolic case, due to the specific form of the Green function $\Gamma$. In particular, for fixed $t$, $\Gamma_t$ has compact support. Throughout this section, we assume that $\sigma(0) = 0$ and that $\Ell_{\sigma} > 0$ (defined in \eqref{eqn:ell}). Hence, the solution to \eqref{eqn:wave} with vanishing initial conditions is identically $0$. We notice that intermittency, as well as position of the peaks and some compact-support related ideas for this case have already been studied in \cite{conus_khoshnevisan}.

The idea borrowed from \cite{foondun_khoshnevisan_2} is to compare $\sup_{x\in\R}u_t(x)$ with the $L^2(\R)$-norm of $u(t,\cdot)$. This comparison will lead to the result since the compact support property of $u_0$ and $v_0$ will lead us to show that $u(t,\cdot)$ has actually compact support. Theorem \ref{thm:sup_compact} below constitutes the main result of this section.

\begin{thm} \label{thm:sup_compact}
Suppose that $\Ell_{\sigma} > 0$ and $\sigma(0)=0$. Assume also that $u_0$ is Hölder-continuous of order $1/2$ and that $u_0$ and $v_0$ are both non-negative functions with compact support included in $[-K,K]$ for some $K > 0$. Let $u$ denote the solution to \eqref{eqn:wave}. Then, $u(t,\cdot) \in L^2(\R)$ a.s. for all $t \geq 0$ and
\begin{equation} \label{eqn:sup_compact}
\Ell_{\sigma} \sqrt{\frac{\kappa}{2}} \leq \limsup_{t \rightarrow \infty} \frac{1}{t} \sup_{x \in \R} \log \E\left[|u(t,x)|^2\right] \leq \limsup_{t \rightarrow \infty} \frac{1}{t} \log \E\left[\sup_{x \in \R} |u(t,x)|^2\right] \leq \Lip_{\sigma} \sqrt{\frac{\kappa}{2}}
\end{equation}
\end{thm}

\begin{rem}
Theorem \ref{thm:sup_compact} implies that the random variable $\sup_{x \in \mathbf{R}} u(t,x)$ is almost surely finite for all finite $t \geq 0$ provided the initial condition have compact support. We are going to show in Section \ref{sec:chaos} that this supremum is almost surely infinite if the initial condition is bounded away from zero. This constitutes an evidence of chaotic behavior of equation \eqref{eqn:wave} since a small perturbation of the initial condition can lead to a drastic change in behavior for the solution.
\end{rem}

Before turning to the proof of Theorem \ref{thm:sup_compact}, we need a few intermediate results.

\begin{prop} \label{prop:sup_compact}
Suppose that $\Ell_{\sigma} > 0$ and $\sigma(0)=0$. Assume that $u_0 \not\equiv 0$ and $v_0$ are non-negative functions in $L^2(\R)$. Let $u$ denote the solution to \eqref{eqn:wave}. Then,
\begin{equation}
\Ell_{\sigma} \sqrt{\frac{\kappa}{2}} \leq \limsup_{t \rightarrow \infty} \frac{1}{t} \log \E\left[\|u(t,\cdot)\|^2_{L^2(\R)}\right] \leq \Lip_{\sigma} \sqrt{\frac{\kappa}{2}}
\end{equation}
\end{prop}

\noindent{\bf Proof.}
This proof follows closely the proof of Theorem 2.1 of \cite{foondun_khoshnevisan}. We refer to it for more details. Since $u_0 \geq 0$, we directly have 
\begin{equation}
\frac{1}{2} \|u_0\|^2_{L^2(\R)} \leq \|U_0(t,\cdot)\|^2_{L^2(\R)} \leq \|u_0\|^2_{L^2(\R)}.
\end{equation}
Moreover, since $v_0 \geq 0$, we have
\begin{eqnarray}
0 \leq \|V_0(t,\cdot)\|^2_{L^2(\R)} & = & \int_{\R} dx \, \left(\int_{-\kappa t}^{\kappa t} dy \, v_0(y+x)\right)^2 \nonumber \\
& \leq & \int_{-\kappa t}^{\kappa t} dy_1 \int_{-\kappa t}^{\kappa t} dy_2 \int_{\R} dx \, v_0(x+y_1)v_0(x+y_2) \nonumber \\
& \leq & 4 \kappa^2 t^2 \, \|v_0\|^2_{L^2(\R)},
\end{eqnarray}
by Hölder's inequality.

Now, from \eqref{eqn:mild}, we have
\begin{eqnarray}
\E\left[\|u(t,\cdot)\|^2_{L^2(\R)}\right] & \geq & \|U_0(t,\cdot)\|^2_{L^2(\R)} + \|V_0(t,\cdot)\|^2_{L^2(\R)} + \Ell^2_{\sigma} \int_{0}^{t} ds \, \E\left[\|u(s,\cdot)\|^2_{L^2(\R)}\right] \|\Gamma_{t-s}\|^2_{L^2(\R)} \nonumber \\
& \geq & \frac{1}{2} \|u_0\|^2_{L^2(\R)} + \Ell^2_{\sigma} \int_{0}^{t} ds \, \E\left[\|u(s,\cdot)\|^2_{L^2(\R)}\right] \|\Gamma_{t-s}\|^2_{L^2(\R)} \label{eqn:compact_1}
\end{eqnarray}
Taking Laplace transforms in \eqref{eqn:compact_1} and setting
\begin{equation}
U(\lambda) := \int_{0}^{\infty} e^{-\lambda t} \E\left[\|u(t,\cdot)\|^2_{L^2(\R)}\right],
\end{equation}
we have, using \eqref{eqn:upsilon},
\begin{equation} \label{eqn:compact_2}
U(\lambda) \geq \frac{\|u_0\|^2_{L^2(\R)}}{2\lambda} + \frac{\kappa \Ell^2_{\sigma}}{2\lambda^2} \, U(\lambda).
\end{equation}
Since $u_0 \not\equiv 0$, the first term on the right-hand side of \eqref{eqn:compact_2} is positive. Hence, this shows that $U(\lambda) = \infty$ provided that $\lambda \leq \Ell_{\sigma} \sqrt{\kappa/2}$ which proves the lower bound in Proposition \ref{prop:sup_compact}.

As for the upper bound, we consider the Picard iteration scheme defining $u$ from \eqref{eqn:mild}, namely
\begin{equation} \label{eqn:picard}
	u_{n+1}(t,x) = U_{0}(t,x) + V_{0}(t,x) + \int_{0}^{t} \int_{\R} \Gamma_{t-s}(y-x) \sigma(u_n(s,y))\, W(ds, dy).
\end{equation}
Similarly as to obtain \eqref{eqn:compact_1}, but using upper bounds, we have
\begin{eqnarray}
\E\left[\|u_{n+1}(t,\cdot)\|^2_{L^2(\R)}\right] & \leq & 2\|U_0(t,\cdot)\|^2_{L^2(\R)} + 2\|V_0(t,\cdot)\|^2_{L^2(\R)} + \Lip^2_{\sigma} \int_{0}^{t} ds \, \E\left[\|u_n(s,\cdot)\|^2_{L^2(\R)}\right] \|\Gamma_{t-s}\|^2_{L^2(\R)} \nonumber \\
& \leq & 2 \|u_0\|^2_{L^2(\R)} + 8\kappa^2 t^2 \, \|v_0\|^2_{L^2(\R)} + \Lip^2_{\sigma} \int_{0}^{t} ds \, \E\left[\|u_n(s,\cdot)\|^2_{L^2(\R)}\right] \|\Gamma_{t-s}\|^2_{L^2(\R)}. \label{eqn:compact_3}
\end{eqnarray}
Now, setting
\begin{equation}
M_n(\lambda) := \sup_{t \geq 0} e^{-\lambda t} \E\left[\|u_n(t,\cdot)\|^2_{L^2(\R)}\right]
\end{equation}
for $n \in \mathbf{N}$, \eqref{eqn:compact_3} and \eqref{eqn:upsilon} lead to
\begin{eqnarray}
M_{n+1}(\lambda) \leq 2 \|u_0\|^2_{L^2(\R)} + \frac{8\kappa^2}{\lambda^2} \, \|v_0\|^2_{L^2(\R)} + \frac{\kappa \Lip^2_{\sigma}}{2 \lambda^2} \, M_n(\lambda).
\end{eqnarray}
Hence, a Gronwall argument shows that $\sup_{n \in \mathbf{N}} M_n(\lambda) < \infty$ provided that $\lambda > \Lip_{\sigma} \sqrt{\kappa/2}$. Taking the limit as $n \rightarrow \infty$ leads to the lower bound in Proposition \ref{prop:sup_compact}. The result is proved.
\hfill $\blacksquare$ \\

Proposition \ref{prop:sup_compact} proves the first claim of Theorem \ref{thm:sup_compact}, namely that the solution $x \mapsto u(t,x)$ is almost surely in $L^2(\R)$ for fixed $t \geq 0$.

We now want to prove Theorem \ref{thm:sup_compact} from Proposition \ref{prop:sup_compact} by showing that $\|u(t,\cdot)\|_{L^2{(\R)}}$ and $\sup_{x\in\R} u(t,x)$ are comparable. We start by a crucial property of the solution $u$. This compares to Lemma 3.3 of \cite{foondun_khoshnevisan}. In the present setting, things are much nicer since the solution $u$ actually has compact support, unlike in the parabolic case, where it only had \emph{essentially} compact support (see \cite{foondun_khoshnevisan} for more details).

\begin{prop} \label{prop:compact}
Under the assumptions of Theorem \ref{thm:sup_compact}, the function $x \mapsto u(t,x)$ has compact support contained in $[-K-\kappa t, K+\kappa t]$. 
\end{prop}

\noindent{\bf Proof.}

As metionned earlier, the solution to \eqref{eqn:wave} with vanishing initial conditions is identically 0. Hence, by Theorem \ref{thm:comp_p}, since $u_0 \geq 0$ and $v_0 \geq 0$, we have $u(t,x) \geq 0$ almost surely for every $t > 0$, $x \in \R$. Now, fix $t > 0$. Since $u_0, v_0$ have support $[-K,K]$, then both $U_0(t,\cdot)$ and $V_0(t,\cdot)$ have support $[-K-\kappa t,K+\kappa t]$. From \eqref{eqn:mild} and since stochastic integrals have zero expectation, we have
\begin{equation}
\E[|u(t,x)|] = \E[u(t,x)] = U_0(t,x) + V_0(t,x),
\end{equation}
and $\E[|u(t,x)|] = 0$ if $|x| > K + \kappa t$. Hence, $u(t,x) = 0$ almost surely for $|x| > K + \kappa t$ and $x \mapsto u(t,x)$ has compact support.
\hfill $\blacksquare$ \\

\begin{rem}
Notice that a consequence of Proposition \ref{prop:compact} is to improve some of the estimates previously obtained in \cite{conus_khoshnevisan}, namely that, for large $t$, $u(t,\cdot)$ does not have large peaks more than a distance $\kappa t$ from the origin.
\end{rem}

In order to be able to prove Theorem \ref{thm:sup_compact}, we need some continuity estimate as regards the solution $u$. We notice that continuity of the solution as been long known (see \cite{walsh, dalang} for instance). We just state the results in the form that we need.

\begin{lemma} \label{lem:cont_1}
Suppose that the initial condition $u_0$ is Hölder-continuous of order $1/2$. Then, for all integers $p \geq 1$ and $\beta > \bar{\gamma}(2p)$, there exists a constant $C_{p,\beta} \in (0,\infty)$ such that, for all $t \geq 0$,
\begin{equation}
\sup_{j \in \mathbf{Z}} \sup_{j \leq x < x' \leq j+1} \left\|\frac{u(t,x)-u(t,x')}{|x-x'|^{1/2}}\right\|_{2p} \leq C_{p,\beta} e^{\beta t/2p}.
\end{equation}
\end{lemma}
 
\noindent{\bf Proof.}
From \eqref{eqn:mild} and the Burkholder-Davis-Gundy inequalities, if $p \geq 1$ is an integer, we have
\begin{eqnarray}
\|u(t,x)-u(t,x')\|_{2p} & \leq & |U_0(t,x)-U_0(t,x')| + |V_0(t,x)-V_0(t,x')| \label{eqn:cont_first} \\
 && \quad + \; 2\sqrt{2p} \, \Lip_{\sigma} \left(\int_{0}^{t} ds \int_{\R} dy \, \|u(s,y)\|^2_{2p} |\Gamma_{t-s}(y-x)-\Gamma_{t-s}(y-x')|^2 \right)^{1/2}. \nonumber
\end{eqnarray}
Now, since $u_0$ is Hölder continuous of order $1/2$, we have
\begin{equation} \label{eqn:cont_u0}
|U_0(t,x)-U_0(t,x')| \leq \const \cdot |x-x'|^{1/2}.
\end{equation}
Moreover, since $v_0$ is bounded, we have
\begin{equation} \label{eqn:cont_v0}
|V_0(t,x)-V_0(t,x')| \leq 2 \overline{v}_0 |x-x'|.
\end{equation}
Finally, a direct calculation leads to
\begin{equation} \label{eqn:cont_gamma}
\int_{\R} dy \, |\Gamma_s(y-x)-\Gamma_s(y-x')|^2 \leq 2 |x-x'|,
\end{equation}
for all $s > 0$. Theorem \ref{thm:inter_up} shows that $\|u\|_{2p,\beta} < \infty$ provided $\beta > \bar{\gamma}(2p)$. As a consequence,
\begin{eqnarray}
\lefteqn{\int_{0}^{t} ds \int_{\R} dy \, \|u(s,y)\|^2_{2p} |\Gamma_{t-s}(y-x)-\Gamma_{t-s}(y-x')|^2} \nonumber \\
\quad & \leq & \|u\|^2_{2p,\beta} \, e^{\beta t/p} \int_{0}^{\infty} ds \, e^{-\beta s/p} \int_{\R} dy \, |\Gamma_{t-s}(y-x)-\Gamma_{t-s}(y-x')|^2 \nonumber \\
\quad & \leq & \|u\|^2_{2p,\beta} \frac{2p}{\beta} \, e^{\beta t/p}  |x-x'|, \label{eqn:cont_global}
\end{eqnarray}
by \eqref{eqn:cont_gamma}. Replacing \eqref{eqn:cont_u0},\eqref{eqn:cont_v0}, and \eqref{eqn:cont_global} in \eqref{eqn:cont_first}, the result follows.
\hfill $\blacksquare$ \\

Similarly as in Lemmas 3.5 and 3.6 of \cite{foondun_khoshnevisan_2}, we can extend this result to all real numbers $p \in (1,2)$ and to a uniform modulus of continuity estimate. We skip the details since they work exactly as in \cite{foondun_khoshnevisan_2}.

\begin{lemma} \label{lem:cont_2}
Suppose the conditions of Lemma \ref{lem:cont_1} are satisfied. Then, for all $p \in (1,2)$ and $\epsilon, \delta \in (0,1)$, there exists a constant $C_{p,\epsilon,\delta} \in (0,\infty)$ such that for all $t \geq 0$,
\begin{equation}
\sup_{j \in \mathbf{Z}} \left\|\sup_{j \leq x < x' \leq j+1} \frac{|u(t,x)-u(t,x')|^2}{|x-x'|^{1-\epsilon}}\right\|_{p} \leq C_{p,\epsilon,\delta} e^{(1+\delta)\lambda_p t},
\end{equation}
where $\lambda_p = (2-p)\bar{\gamma}(2) + (p-1)\bar{\gamma}(4)$.
\end{lemma}

We are now ready to prove Theorem \ref{thm:sup_compact}. This is similar to the proof of Theorem 1.1 in \cite{foondun_khoshnevisan_2}, but this case is easier due to Proposition \ref{prop:compact}. \\

\noindent{\bf Proof of Theorem \ref{thm:sup_compact}:}
We already proved that $u(t,\cdot) \in L^2(\R)$ with Proposition \ref{prop:sup_compact}. It remains to prove \eqref{eqn:sup_compact}

The lower bound is a direct consequence of Propositions \ref{prop:sup_compact} and \ref{prop:compact}. Indeed, from Proposition \ref{prop:sup_compact}, we have
\begin{eqnarray}
\exp\left(\left[\Ell_{\sigma} \sqrt{\frac{\kappa}{2}} + o(1)\right] t\right) & \leq & \E\left[\int_{\R} |u(t,x)|^2 \, dx\right] \nonumber \\
& = & \E\left[\int_{-K-\kappa t}^{K + \kappa t} |u(t,x)|^2 \, dx\right] \nonumber \\
& \leq & 2(K + \kappa t) \sup_{x \in \R} \E[|u(t,x)|^2].
\end{eqnarray}
The first inequality in \eqref{eqn:sup_compact} follows.

As for the upper bound, for all $p \in (1,2)$, $\epsilon \in (0,1)$, $j \in \mathbf{Z}$ and $t \geq 0$, we have
\begin{eqnarray}
\sup_{j \leq x \leq j+1} |u(t,x)|^{2p} & \leq & 2^{2p-1} \left(|u(t,j)|^{2p} + \sup_{j \leq x \leq j+1} |u(t,x)-u(t,j)|^{2p} \right) \nonumber \\
& \leq & 2^{2p-1}\left(|u(t,j)|^{2p} + \Omega_j^p\right).
\end{eqnarray}
where
\begin{equation}
\Omega_j^p := \sup_{j \leq x \leq x' \leq j+1} \frac{|u(t,x)-u(t,x')|^2}{|x-x'|^{1-\epsilon}}.
\end{equation}
It follows that
\begin{equation}
\E\left[\sup_{j \leq x \leq j+1} |u(t,x)|^{2p}\right] \leq 2^{2p-1} \left(\E[|u(t,j)|^{2p}] + \E[\Omega_j^p]\right).
\end{equation}
Now, Lemma \ref{lem:cont_2} implies that $\E[\Omega_j^p] \leq C_{p,\epsilon,\delta} e^{p(1+\delta)\lambda_p t}$. Moreover, $\E[u(t,j)^{2p}] = 0$ for $|j| > K+\kappa t$, and by Theorem \ref{thm:inter_up}, $\E[u(t,j)^{2p}] \leq \const \cdot e^{(\bar{\gamma}(2p)+o(1))t}$ for $|j| \leq K + \kappa t$. It follows that
\begin{eqnarray}
\E\left[\sup_{x \in \R} |u(t,x)|^{2p}\right] & \leq & \E\left[\sup_{|x| \leq \lceil K+\kappa t\rceil} |u(t,x)|^{2p}\right] \nonumber \\
& \leq & \const \cdot \lceil K+\kappa t\rceil \left(e^{(\bar{\gamma}(2p)+o(1))t} + C_{p,\epsilon,\delta} e^{p(1+\delta)\lambda_p t}\right).
\end{eqnarray}
Hence,
\begin{equation}
\limsup_{t \rightarrow \infty} \frac{1}{t} \, \log \E\left[\sup_{x \in \R} |u(t,x)|^{2p}\right] \leq \max\{p(1+\delta)\lambda_p ; \bar{\gamma}(2p)\}.
\end{equation}
Taking $\delta \rightarrow 0$, then using Jensen's inequality before taking $p \rightarrow 1$ leads to
\begin{equation}
\limsup_{t \rightarrow \infty} \frac{1}{t} \, \log \E\left[\sup_{x \in \R} |u(t,x)|^{2}\right] \leq \bar{\gamma}(2) \leq \Lip_{\sigma} \sqrt{\frac{\kappa}{2}},
\end{equation}
by Theorem \ref{thm:inter_up} and Remark \ref{rem:inter_up}. This proves the last inequality in \eqref{eqn:sup_compact} and concludes the proof.
\hfill $\blacksquare$


\section{Moment and tail probability estimates} \label{sec:moment_tail}

In this section, we will present technical estimates, first on moments, which will then lead to estimates on tail probabilities. These results will be used to prove the main results of this paper in Section \ref{sec:chaos}. The results of this section are comparable to the ones of \cite{CJK} for the parabolic equations. Hence, we mainly concentrate on the differences in the proofs below and refer to \cite{CJK} for more details. 

Throughout Section \ref{sec:moment_tail}, we will consider constant initial conditions $u_0$ and $v_0$. The reasons for this restriction are outlined in Section \ref{sec:chaos}. Without loss of generality, we will assume that $u_0 \equiv 1$. We keep track of $v_0$, since it makes a difference below. Hence, \eqref{eqn:mild} becomes
\begin{equation} \label{eqn:mild_c} 
u(t,x) = 1 + v_0 \kappa t + \int_{0}^{t} \int_{\R} \Gamma_{t-s}(y-x) \sigma(u(s,y)) W(ds,dy),
\end{equation}
for $t \geq 0$, $x \in \R$. We notice that the results of Dalang \cite{dalang} imply that the law of $u_t(x)$ is independent of $x$ when the initial conditions are constant.

Let us start our presentation by proving a general upper bound for moments. 

\begin{prop} \label{prop:mom_ub}
Let $u$ denote the solution to \eqref{eqn:wave} with initial conditions identically one. Fix $T > 0$ and let $a := T \Lip_{\sigma} \sqrt{\kappa}$, then, for every $p \geq 1$, there exists a constant $C > 0$, independent of $p$ such that
\begin{equation}
\sup_{0 \leq t \leq T} \sup_{x \in \R} \E[|u(t,x)|^p] \leq C^{p} e^{a p^{3/2}}.
\end{equation}
\end{prop}

\noindent{\bf Proof.}
This is a direct consequence of the proof of Theorem \ref{thm:inter_up}. Indeed, we proved there that $\|u\|_{p,\beta} < \infty$ provided $\beta > p^{3/2} \Lip_{\sigma} \sqrt{\kappa/2}$. The result follows.
\hfill $\blacksquare$ \\

Now, we would like to turn this moment estimate into an upper bound on the tail of the distribution of $u_t(x)$. Corollary \ref{coro:prob_ub} is a direct consequence of Proposition \ref{prop:mom_ub} and Lemma 3.4 of \cite{CJK}. We will skip the details and refer to \cite{CJK} for more on this machinery.

\begin{coro} \label{coro:prob_ub}
Fix $T > 0$. Then, for all $\alpha < \frac{4}{27} (T^2 (\Lip_{\sigma} \vee 1)^2 \kappa)^{-1}$,
\begin{equation}
\sup_{0 \leq t \leq T} \sup_{x \in \R} \E\left[\exp\left(\alpha (\log_+(u(t,x))^3 \right)\right] < \infty,
\end{equation}
where $\log_+(x) = \log(x \vee e)$. As a consequence,
\begin{equation}
\limsup_{\lambda \rightarrow \infty} \frac{1}{(\log \lambda)^3} \sup_{0 \leq t \leq T} \sup_{x \in \R} \log \P\{u(t,x) > \lambda\} \leq - \frac{4}{27 \, T^2 (\Lip_{\sigma} \vee 1)^2 \kappa}.
\end{equation}
\end{coro}

In order to control the moments and, hence, the tail probabilities, we will need a lower bound equivalent result to Proposition \ref{prop:mom_ub} and Corollary \ref{coro:prob_ub}. This bound will depend on the behavior of the non-linearity $\sigma$. Let us start with the easiest case, where $\sigma$ is bounded away from 0.

\begin{prop} \label{prop:mom_lb_bdd}
Let $\epsilon_0 := \inf_{z \in \R} \sigma(z) > 0$. Then, for all $t > 0$ and all $p \geq 1$,
\begin{equation} \label{eqn:mom_lb_bdd}
\inf_{x \in \R} \E[|u(t,x)|^{2p}] \geq (\sqrt{2} + o(1)) (\mu_t p)^p \qquad (\mbox{as } p \rightarrow \infty),
\end{equation}
where the $o(1)$ term only depends on $p$ and
\begin{equation} \label{eqn:mu_t}
\mu_t := \frac{\epsilon_0^2 \kappa t^2}{2e}.
\end{equation}
\end{prop}

\noindent{\bf Proof.}
Since the law of $u(t,x)$ doesn't depend on $x$, the $\inf$ in \eqref{eqn:mom_lb_bdd} is not needed and, without loss of generality, we can assume $x=0$. For a fixed $t > 0$, we first notice that $u(t,0) = 1 + v_0 \kappa t + M_t$, where $(M_{\tau})_{0 \leq \tau \leq t}$ is the martingale defined by
\begin{equation}
M_{\tau} := \int_{0}^{\tau} \int_{\R} \Gamma_{t-s}(y) \sigma(u(s,y)) W(ds,dy).
\end{equation}
Now, the quadratic variation of $M$ is given by
\begin{equation} \label{eqn:qv}
\langle M\rangle_{\tau} = \int_{0}^{\tau} ds \int_{\R} dy \, \Gamma^2_{t-s}(y) \sigma^2(u(s,y)).
\end{equation}
An iterative use of Itô's formula, \eqref{eqn:qv} and the fact that $\sigma(z) \geq \epsilon_0$ for all $z \in \R$ leads to
\begin{equation} \label{eqn:mom_formula}
\E[M_t^{2p}] \geq \sum_{k=0}^{p-1} \left(\prod_{j=0}^{k} \binom{2p-2j}{2}\right) \epsilon_0^{2(k+1)} \int_{0}^{t} \nu(t,ds_1) \int_{0}^{s_1} \nu(s_1,ds_2) \cdots \int_{0}^{s_k} \nu(s_k,ds_{k+1}),
\end{equation}
where the measure $\nu(t,ds)$ is given by
\begin{equation}
\nu(t,ds) := \ind_{[0,t]}(s) \, \|\Gamma_{t-s}\|^2_{L^2(\R)} \, ds.
\end{equation}
We refer to the proof of Lemma 3.6 of \cite{CJK} for a detailed analogue in the parabolic case. We also refer to \cite{conus_dalang} for similar moment computations in the hyperbolic case. The right-hand side of \eqref{eqn:mom_formula} is the exact expression for the $p^{\mbox{\scriptsize th}}$ moment of $u$ if $\sigma$ were identically $\epsilon_0$. Hence, we have shown that
\begin{equation} \label{eqn:bound}
\E[|u(t,0)|^{2p}] \geq \E[M_t^{2p}] \geq \E[N_t^{2p}],
\end{equation}
where
\begin{equation}
N_t := \epsilon_0 \, \int_{0}^{t} \int_{\R} \Gamma_{t-s}(y) W(ds,dy).
\end{equation}
Since $N_t$ is Gaussian and has second moment given by
\begin{equation}
\E[N_t^2] = \epsilon_0^2 \int_{0}^{t} ds \int_{\R} dy \, \Gamma^2_{s}(y) = \frac{\epsilon_0^2 \kappa t^2}{4},
\end{equation}
we have
\begin{equation} \label{eqn:mom_final}
\E[N_t^{2p}] = \frac{(2p)!}{2^p \, p!} \E[N_t^2]^p = \frac{(2p)!}{2^p \, p!} \left(\frac{\epsilon_0^2 \kappa t^2}{4}\right)^p.
\end{equation}
Stirling's formula, \eqref{eqn:bound} and \eqref{eqn:mom_final} prove the result.
\hfill $\blacksquare$ \\

Now, we would like to turn the lower bound on moments of Proposition \ref{prop:mom_lb_bdd} into lower bounds on the tail probabilities. This uses the so-called Paley-Zygmund inequality similarly as in the proof of Proposition 3.7 of \cite{CJK}.

\begin{prop} \label{prop:pr_lb_bdd}
Let $\inf_{z \in \R} \sigma(z) = \epsilon_0 > 0$. Then, there exists a constant $C \in (0,\infty)$ such that, for all $t > 0$,
\begin{equation}
\liminf_{\lambda \rightarrow \infty} \frac{1}{\lambda^3} \inf_{x \in \R} \log \P\{|u(t,x)| \geq \lambda\} \geq - C \, \frac{(\Lip_{\sigma} \vee 1)}{\epsilon_0^3 t^2 \kappa}.
\end{equation}
\end{prop}

\noindent{\bf Proof.}
The Paley-Zygmund inequality (a derivation is proposed in the proof of Proposition 3.7 in \cite{CJK}) states that
\begin{eqnarray}
\P\left\{|u(t,x)| \geq \frac{1}{2} \|u(t,x)\|_{2p}\right\} & \geq & \frac{\E\left[|u(t,x)|^{2p}\right]^2}{4 \E\left[|u(t,x)|^{4p}\right]} \\
& \geq & \exp\left(-8 t (\Lip_{\sigma} \vee 1) \kappa^{1/2} \, p^{3/2}\right),
\end{eqnarray}
by the bounds of Propositions \ref{prop:mom_ub} and \ref{prop:mom_lb_bdd}. Moreover, Proposition \ref{prop:mom_lb_bdd} shows that
\begin{equation}
\|u(t,x)\|_{2p} \geq (1 + o(1))(\mu_t p)^{1/2},
\end{equation}
as $p \rightarrow \infty$, where $\mu_t$ is given by \eqref{eqn:mu_t}. This implies that
\begin{equation}
\P\left\{|u(t,x)| \geq \frac{1}{2} (\mu_t p)^{1/2}\right\} \geq \exp\left(-8 t (\Lip_{\sigma} \vee 1) \kappa^{1/2} \, p^{3/2}\right)
\end{equation}
as $p \rightarrow \infty$. Considering $\lambda := (\mu_t p)^{1/2}$, the result follows.
\hfill $\blacksquare$ \\

The results above were obtained under the condition that $\sigma$ was bounded away from $0$. In the case were $\sigma(z)$ decreases to $0$ not too fast as $|z| \rightarrow \infty$, we can still obtain similar lower bounds. Namely, consider $\sigma$ to satisfy
\begin{equation} \label{eqn:sigma_log}
\lim_{|z| \rightarrow \infty} \sigma(z) \log(|z|)^{\frac{1}{3}-\gamma} = \infty,
\end{equation}
for some $\gamma \in \left(0,\frac{1}{3}\right)$. Then, we obtain the lower bound given in Proposition \ref{prop:pr_lb_log} below. The proof follows exactly the arguments of Proposition 3.8 in \cite{CJK} (using Proposition \ref{prop:pr_lb_bdd} instead of the parabolic equivalent) and we skip the details.

\begin{prop} \label{prop:pr_lb_log}
Assume $\sigma$ satisfies \eqref{eqn:sigma_log} for some $\gamma \in \left(0,\frac{1}{3}\right)$. Then, there exists a constant $C \in (0,\infty)$ depending only on $\gamma$, such that for all $t > 0$,
\begin{equation}
\liminf_{\lambda \rightarrow \infty} \frac{1}{\lambda^{1/\gamma}} \inf_{x \in \R} \log \P\{|u(t,x)| \geq \lambda\} \geq - C \left(\frac{(\Lip_{\sigma} \vee 1)}{t^2 \kappa}\right)^{1/3\gamma}.
\end{equation}
\end{prop}

Summarizing the results obtained so far, we obtain Corollary \ref{coro:bounds} below. We write $f(x) \succsim g(x)$ as $x \rightarrow \infty$ instead of ``there exists a deterministic constant $C$ such that $\liminf_{x \rightarrow \infty} f(x)/g(x) \geq C$''.

\begin{coro} \label{coro:bounds}
Let $u$ denote the solution to \eqref{eqn:wave} with initial conditions identically one. If $\inf_{z \in \R} \sigma(z) = \epsilon_0 > 0$, then for all $t > 0$,
\begin{equation}
- \frac{\lambda^3}{\kappa} \precsim \log \P\{|u(t,x)| \geq \lambda\} \precsim - \frac{(\log\lambda)^3}{\kappa}, \qquad \mbox{ as } \lambda \rightarrow \infty.
\end{equation}
If $\sigma$ satisfies \eqref{eqn:sigma_log} for some $\gamma \in \left(0,\frac{1}{3}\right)$, then for all $t > 0$,
\begin{equation}
- \frac{\lambda^{1/\gamma}}{\kappa^{1/3\gamma}} \precsim \log \P\{|u(t,x)| \geq \lambda\} \precsim - \frac{(\log\lambda)^3}{\kappa}, \qquad \mbox{ as } \lambda \rightarrow \infty.
\end{equation}
The inequalities above hold uniformly for all $x \in \R$ and the constants behind $\precsim$ do not depend on $\kappa$. 
\end{coro}

The upper bounds obtained in Proposition \ref{prop:mom_ub} and Corollary \ref{coro:prob_ub} were pretty general and only assumed that $\sigma$ was a Lipschitz function. Now, if we assume that $\sigma$ is bounded above (as well as bounded away from $0$), these are far from optimal. Actually, we can show that the lower bound of Proposition \ref{prop:mom_lb_bdd} is sharp in that case.

\begin{prop} \label{prop:mom_ub_bdd}
Let $S_0 := \sup_{z \in \R} \sigma(z) < \infty$. Then, for all $t > 0$ and all integers $p \geq 1$,
\begin{equation} \label{eqn:mom_ub_bdd}
\sup_{x \in \R} \E[|u(t,x)|^{2p}] \leq (2\sqrt{2} + o(1)) (\tilde{\mu}_t p)^p \qquad (\mbox{as } p \rightarrow \infty),
\end{equation}
where the $o(1)$ term only depends on $p$ and
\begin{equation} \label{eqn:mu_t_tilde}
\tilde{\mu}_t := \max\left(\frac{2 S_0^2 \kappa t^2}{e}, 4 v_0^2 \kappa^2 t^2\right).
\end{equation}
\end{prop}

\noindent{\bf Proof.}
We follow an argument similar to the one in the proof of Proposition \ref{prop:mom_lb_bdd}. We consider the same martingale $(M_{\tau})_{0\leq\tau\leq t}$ as in Proposition \ref{prop:mom_lb_bdd}. The exact same argument, but reversing the inequalities and using $\sigma(z) \leq S_0$ for all $z \in \R$, shows that
\begin{equation}
\E[|u(t,0)|^{2p}] \leq 2^{2p}(1+v_0 \kappa t)^{2p} + 2^{2p} \E[M_t^{2p}] \leq 2^{2p}(1+v_0 \kappa t)^{2p} + 2^{2p} \E[N_t^{2p}],
\end{equation}
where $(N_t)_{t \geq 0}$ is defined by
\begin{equation}
N_t := S_0 \, \int_{0}^{t} \int_{\R} \Gamma_{t-s}(y) W(ds,dy).
\end{equation}
Similar computations as in Proposition \ref{prop:mom_lb_bdd} prove the result.
\hfill $\blacksquare$ \\

We can now turn this bound into estimates on the probability tail.

\begin{prop} \label{prop:pr_bdd}
Let $u$ be the solution to \eqref{eqn:wave} with $0 < \epsilon_0 := \inf_{z \in \R} \sigma(z) \leq \sup_{z \in \R} \sigma(z) := S_0 < \infty$. Then, for all $t > 0$, there exists constants $C > c > 0$ such that, simultaneously for all $\lambda$ large enough and $x \in \R$,
\begin{equation}
c \exp\left(-C\frac{\lambda^2}{\kappa}\right) \leq \P\{|u(t,x)| \geq \lambda\} \leq C \exp\left(-c\frac{\lambda^2}{\max\{\kappa,v_0^2 \kappa^2\}}\right).
\end{equation}
\end{prop}

\begin{rem}
Notice that if $v_0 \equiv 0$, then the behavior in $\kappa$ of both the upper and lower bound agree.
\end{rem}

\noindent{\bf Proof.}
The lower bound is obtained in the exact same way as in Proposition \ref{prop:pr_lb_bdd} using the Paley-Zigmund inequality, but replacing the moment upper bound of Proposition \ref{prop:mom_ub} by the one obtained above in Proposition \ref{prop:mom_ub_bdd}.

As for the lower bound, from Proposition \ref{prop:mom_ub_bdd}, we have that for all integers $p \geq 0$, $\sup_{x \in \R} \E[|u(t,x)|^{2p}] \leq (A\max\{\kappa,v_0^2 \kappa^2\})^m \, m!$, for some constant $A \in (0,\infty)$. This implies that
\begin{equation}
\sup_{x \in \R} \E\left[\exp\left(\xi |u(t,x)|^2\right)\right] \leq \sum_{p=0}^{\infty} (\xi A \max\{\kappa,v_0^2 \kappa^2\})^m = \frac{1}{1-\xi A \max\{\kappa,v_0^2 \kappa^2\}} < \infty,
\end{equation}
for $\xi < (A\max\{\kappa,v_0^2 \kappa^2\})^{-1}$. Then, for such a $\xi$, Chebychev's inequality implies 
\begin{equation}
\P\{|u(t,x)| > \lambda\} \leq \frac{\exp(-\xi \lambda^2)}{1-\xi A \max\{\kappa,v_0^2 \kappa^2\}}.
\end{equation} 
Since this is valid for all $x \in \R$, we choose $\xi = \const \cdot (\max\{\kappa,v_0^2 \kappa^2\})^{-1}$ to obtain the result.
\hfill $\blacksquare$


\section{Localization} \label{sec:loc}

One of the main argument that will lead to estimates on the supremum of the solution to \eqref{eqn:wave} in Section \ref{sec:chaos} is the so-called \emph{localization} property. This property essentially states that if $x_1$ and $x_2$ are chosen sufficiently far apart, then $u(t,x_1)$ and $u(t,x_2)$ are approximately independent. This idea was already used in \cite{CJK} for the study of the stochastic heat equation. In \cite{CJK}, a precise estimate based on the exponential decrease property of the heat kernel was needed together with the independent increment property of the space-time white noise. In our case, this matter is made much easier by the compact support property of $\Gamma$. Below, we will state the localization results for the hyperbolic equation \eqref{eqn:wave}. Similarly as in Section \ref{sec:moment_tail}, we will assume that the initial conditions are identically constant and that $u_0 \equiv 1$; we refer to Section \eqref{sec:chaos} for the reasons of this restriction.

We remind that when the initial condition are identically constant, the solution to \eqref{eqn:wave} satisfies the mild form given by
\begin{equation} \label{eqn:mild_c2} 
u(t,x) = 1 + v_0 \kappa t + \int_{0}^{t} \int_{\R} \Gamma_{t-s}(y-x) \sigma(u(s,y)) W(ds,dy),
\end{equation}
for $t \geq 0$, $x \in \R$.

Now, by \eqref{eqn:green_wave}, $\Gamma_{t-s}(z) = 0$, provided $|z| > \kappa(t-s)$. Hence, the support of the space integral in \eqref{eqn:mild_c2} is given by $\{y \in \R : |y-x|\leq \kappa(t-s)\}$, which is contained in $[x-\kappa t, x + \kappa t]$. As a consequence, we have
\begin{equation} \label{eqn:mild_kappa} 
u(t,x) = 1 + v_0 \kappa t + \int_{0}^{t} \int_{x-\kappa t}^{x + \kappa t} \Gamma_{t-s}(y-x) \sigma(u(s,y)) W(ds,dy),
\end{equation}
for all $t \geq 0$, $x \in \R$. Now, for all $n \in \mathbf{N}$, let $\{u_n(t,x); t \geq 0, x \in \R\}$ be the $n$-th step Picard approximation to $u$. Namely, we have $u_0 \equiv 0$ and, for $n \geq 1$, $t \geq 0$ and $x \in \R$,
\begin{equation} \label{eqn:mild_Picard} 
u_n(t,x) = 1 + v_0 \kappa t + \int_{0}^{t} \int_{x-\kappa t}^{x + \kappa t} \Gamma_{t-s}(y-x) \sigma(u_{n-1}(s,y)) W(ds,dy).
\end{equation}

\begin{prop} \label{prop:Picard}
Let $u$ be the solution to \eqref{eqn:wave} with initial conditions identically $1$ and $u_n$ as defined above. Then, for all $t \geq 0$, we have
\begin{equation}
\sup_{x \in \R} \E\left[|u(t,x)-u_n(t,x)|^p\right] \leq C^p e^{a p^{3/2} t} e^{-np},
\end{equation}
where $C$ and $a$ are the constants of Proposition \ref{prop:mom_ub}. 
\end{prop}

\noindent{\bf Proof.}
Proposition \ref{prop:young}, together with \eqref{eqn:mild_c2} and \eqref{eqn:mild_Picard}, show that
\begin{equation}
\|u-u_n\|_{p,\beta} \leq \const \cdot \sqrt{p} \, \Upsilon\left(\frac{2\beta}{p}\right)^{1/2} \|u-u_{n-1}\|_{p,\beta}.
\end{equation}
Hence, using \eqref{eqn:upsilon}, if we choose $\beta = D p^{3/2}$ for a sufficiently large constant $D$ (compare with \eqref{eqn:beta_large}), we can show that
\begin{equation} \label{eqn:e-1}
\|u-u_n\|_{p,\beta} \leq e^{-1} \|u-u_{n-1}\|_{p,\beta}.
\end{equation}
The result follows from \eqref{eqn:e-1} and Theorem \ref{thm:inter_up}.
\hfill $\blacksquare$ \\

Now, we use the fact that the $n$-th Picard approximation $u_n(t,x)$ only depends on the noise $W(s,y)$ for $s \in [0,t]$ and $y \in [x-n\kappa t, x + n\kappa t]$. (We can easily prove this by induction, see \cite[Lemma 4.4 and Appendix A]{CJK}.) This and the properties of stochastic integrals with respect to space-time white noise lead to Proposition \ref{prop:indep} below. The proof follows Lemma 4.4 of \cite{CJK} and we skip the details.

\begin{prop} \label{prop:indep}
Let $t > 0$ and choose $n \in \mathbf{N}$. Now, let $(x_i)_{i \in \mathbf{N}}$ be a sequence such that $|x_i-x_j| \geq 2 n \kappa t$, whenever $i \neq j$. Then $\{u_n(t,x_i)\}_{i \in \mathbf{N}}$ is a collection of i.i.d. random variables.
\end{prop}


\section{Chaotic behavior} \label{sec:chaos}

We are now ready to state and prove the main results of this paper, Theorems \ref{thm:chaos_main} and \ref{thm:chaos_bdd} below. We will use the results from Sections \ref{sec:moment_tail} and \ref{sec:loc} above. The proofs follows similar ideas as the proof of Theorems 1.1 and 1.2 of \cite{CJK}. Nevertheless, we will still give some details for the sake of completeness. We write $f(x) \succsim g(x)$ as $x \rightarrow \infty$ instead of ``there exists a deterministic constant $C$ such that $\liminf_{x \rightarrow \infty} f(x)/g(x) \geq C$''.

\begin{thm} \label{thm:chaos_main}
Let $u$ be the solution to \eqref{eqn:wave} with initial conditions satisfying
\begin{equation}
\inf_{x \in \R} u_0(x) > 0 \qquad \mbox{ and } \qquad v_0(x) \geq 0, \mbox{ for all } x \in \R.
\end{equation}
Then, the following hold:
\begin{itemize}
\item[1.] If $\inf_{z \in \R} \sigma(z) = \epsilon_0 > 0$ and $t > 0$, then
\begin{equation}
\sup_{x \in [-R,R]} u(t,x) \succsim \kappa^{1/3} (\log R)^{1/3} \qquad \mbox{ a.s. as } R \rightarrow \infty.
\end{equation}
\item[2.] If there exists $\gamma \in \left(0,\frac{1}{3}\right)$ such that
\begin{equation} \label{eqn:sigma_log_thm}
\lim_{|z| \rightarrow \infty} \sigma(z) \log(|z|)^{\frac{1}{3}-\gamma} = \infty,
\end{equation}
then, for all $t > 0$,
\end{itemize}
\begin{equation}
\sup_{x \in [-R,R]} u(t,x) \succsim \kappa^{1/3} (\log R)^{\gamma} \qquad \mbox{ a.s. as } R \rightarrow \infty.
\end{equation}
\end{thm}

Similarly as in the parabolic case (see \cite{CJK}), this result establishes a rate of blow-up of $(\log R)^{1/3}$ which is independent of $\sigma$ and the initial conditions in the first part. We notice that this rate is actually different from the one obtained in the parabolic case (namely, $(\log R)^{1/6}$). We would like to point out the dependence with respect to $\kappa$ which is drastically different from the one obtained in \cite{CJK}. First of all, the supremum in Theorem \ref{thm:chaos_main} gets smaller as $\kappa$ goes to $0$, unlike in the parabolic case. Such a difference was already noticed in \cite{conus_khoshnevisan} about the position of the peaks in the case of a compact support initial data. Moreover, the relation between the powers of $\log R$ and $\kappa$ is of the same exponential order ($(\log R)^{1/3}$ and $\kappa^{1/3}$) unlike the parabolic case where we had $(\log R)^{1/6}$ and $\kappa^{1/12}$. As $\kappa$ corresponds in some sense to $1/t$, this suggests the asymptotic space-time scaling behavior of $x \sim \sqrt{t}$ for the parabolic equation and $x \sim t$ in our present hyperbolic case.

Before we turn to the proof, we would like to mention that it is sufficient to prove the result in the case where $u_0$ and $v_0$ are constant. Indeed, since $0 < \underline{u}_0 \leq u_0(x) \leq \overline{u}_0 < \infty$ and $0 \leq v_0(x) \leq \overline{v}_0 < \infty$, the comparison principle developped in Section \ref{sec:comp} (Theorem \ref{thm:comp_p}) will prove the result as soon as it is proved in the case where $u_0 \simeq \underline{u}_0$ and $v_0 \simeq 0$. \\

\noindent{\bf Proof.}
As mentioned above, we only consider the case where the initial conditions are constant. Hence, the results of Section \ref{sec:moment_tail} apply. We will only present in detail the proof of the second part of Theorem \ref{thm:chaos_main}. For the first case, it suffices to take $\gamma = 1/3$ in the argument below. Fix integers $n,N > 0$ and let $(x_i)_{i=1}^{N}$ be a sequence of points as in Proposition \ref{prop:indep}. Then, by Proposition \ref{prop:indep}, $(u_n(t,x_i))_{i=1}^{N}$ is a sequence of independent random variables.
Let $\lambda > 0$, we have
\begin{eqnarray}
\P\left\{\max_{1 \leq j \leq N} |u(t,x_j)| < \lambda\right\} \leq \P\left\{\max_{1 \leq j \leq N} |u_n(t,x_j)| < 2\lambda\right\} + \P\left\{\max_{1 \leq j \leq N} |u(t,x_j)-u_n(t,x_j)| > \lambda\right\}.
\end{eqnarray}
Now, we can apply Proposition \ref{prop:pr_lb_log} (which easily generalizes to $u_n$) and the independence of the random variables to obtain
\begin{equation}
\P\left\{\max_{1 \leq j \leq N} |u_n(t,x_j)| < 2\lambda\right\} \leq \left(1 - c_1 e^{-c_2 (2\lambda)^{1/\gamma}}\right)^N,
\end{equation}
for some constants $c_1$ and $c_2$. Moreover, Chebychev's inequality together with Proposition \ref{prop:Picard} shows that
\begin{equation}
\P\left\{\max_{1 \leq j \leq N} |u(t,x_j)-u_n(t,x_j)| > \lambda\right\} \leq N C^p e^{a p^{3/2} t} e^{-np} \lambda^{-p}.
\end{equation}
Hence, 
\begin{equation} \label{eqn:main_1}
\P\left\{\max_{1 \leq j \leq N} |u(t,x_j)| < \lambda\right\} \leq \left(1 - c_1 e^{-c_2 (2\lambda)^{1/\gamma}}\right)^N + N C^p e^{a p^{3/2} t} e^{-np} \lambda^{-p}.
\end{equation}
Now, we choose the parameters judiciously: we take $\lambda:=p$, $N:=p e^{c_2 p^{1/\gamma}}$, $n=\varrho p^{(1-\gamma)/\gamma}$, for some constant $\varrho > 2^{1/\gamma} c_2$.
As a consequence, \eqref{eqn:main_1} becomes
\begin{eqnarray}
\P\left\{\max_{1 \leq j \leq N} |u(t,x_j)| < p\right\} & \leq & e^{-c_1 p} + \exp\left(c_2(2p)^{1/\gamma} + \log(p) + at p^{3/2} -\varrho p^{1/\gamma} - p\log(p)\right) \nonumber \\
& \leq & 2 e^{-c_1 p},
\end{eqnarray}
since $1/\gamma > 3$. Now, we can choose $x_i = 2 i \kappa t n$, which together with a symmetry argument leads to
\begin{equation} \label{eqn:main_2}
\P\left\{\sup_{|x|\leq 2 N \kappa t n} |u(t,x)| < p\right\} \leq 2 e^{-c_1 p}.
\end{equation}
Now, as $p \rightarrow \infty$,
\begin{equation}
2N\kappa t n  = O(e^{c_2 p^{1/\gamma}}).
\end{equation}
The Borel-Cantelli lemma, together with a monotonicity argument shows that 
\begin{equation}
\sup_{|x| < R} u(t,x) \geq \const \cdot \left(\log(R)/c_2\right)^{\gamma}.
\end{equation} 
Now, by Proposition \ref{prop:pr_lb_log}, $c_2 = \const \cdot \kappa^{-1/3\gamma}$. The result follows.
\hfill $\blacksquare$ \\

Now, we would like to study the case where $\sigma$ is bounded away from $0$ (as above), but also bounded above. We would like to show that the behavior of the solution $u$ in that case is essentially similar to the case where $\sigma$ is identically constant, in which $u(t,x)$ is a Gaussian process. 

\begin{thm} \label{thm:chaos_bdd}
Let $u$ be the solution to \eqref{eqn:wave} with initial conditions satisfying
\begin{equation}
\inf_{x \in \R} u_0(x) > 0 \qquad \mbox{ and } \qquad v_0(x) \geq 0, \mbox{ for all } x \in \R.
\end{equation}
Moreover, assume $0 < \inf_{z \in \R} \sigma(z) \leq \sup_{z \in \R} \sigma(z) < \infty$. Then, for all $t > 0$,
\begin{equation}
\kappa^{1/2} (\log R)^{1/2} \, \precsim \, \sup_{x \in [-R,R]} u(t,x) \, \precsim \, \max\{\kappa^{1/2},\overline{v}_0 \kappa\} \cdot (\log R)^{1/2}, \qquad \mbox{a.s. as } R \rightarrow \infty,
\end{equation}
where the constants behind $\precsim$ do not depend on $\kappa \geq \kappa_0$ for some appropriate constant $\kappa_0$.
\end{thm} 

Before we turn to the proof of Theorem \ref{thm:chaos_bdd}, we notice that as in Theorem \ref{thm:chaos_main}, we only need to prove the result for constant $u_0$ and $v_0$. Indeed, since $0 < \underline{u}_0 \leq u_0(x) \leq \overline{u}_0 < \infty$ and $0 \leq v_0(x) \leq \overline{v}_0 < \infty$, the comparison principle (Theorem \ref{thm:comp_p}) will prove the result as soon as it is proved in the case where 1. $u_0 \simeq \underline{u}_0$ and $v_0 \simeq 0$ for the lower bound and 2. $u_0 \simeq \overline{u}_0$ and $v_0 \simeq \overline{v}_0$ for the upper bound.

We first need a spatial continuity estimate for the solution $u$, simlarly as Lemma 6.1 in \cite{CJK}.

\begin{lemma}
Let $u$ be the solution to \eqref{eqn:wave} with $0 < \epsilon_0 := \inf_{z \in \R} \sigma(z) \leq \sup_{z \in \R} \sigma(z) := S_0 < \infty$. Then, for every $t > 0$, there exists a constant $A \in (0,\infty)$ such that, for all $p \geq 2$,
\begin{equation}
\sup_{x\neq x'} \frac{\E\left[|u(t,x)-u(t,x')|^2p\right]}{|x-x'|^p} \leq (A p)^p.
\end{equation}
\end{lemma}

\noindent{\bf Proof.}
We proceed as in the proof of Lemma 6.1 of \cite{CJK}. Fix $x, x' \in \R$ and let $(M_{\tau})_{0\leq \tau \leq t}$ be the martingale defined by
\begin{equation}
M_{\tau} := \int_{0}^{\tau} \int_{\R} (\Gamma_{t-s}(y-x)-\Gamma_{t-s}(y-x')) \sigma(u(s,y)) W(ds,dy).
\end{equation}
Then, we clearly have
\begin{eqnarray}
\langle M \rangle_{\tau} & \leq & S_0^2 \int_{0}^{\tau} ds \int_{\R} dy \, |\Gamma_{t-s}(y-x)-\Gamma_{t-s}(y-x')|^2 \nonumber \\
& \leq & 2 \tau S_0^2 |x-x'|,  
\end{eqnarray}
by \eqref{eqn:cont_gamma}. The Burkholder-Davis-Gundy inequality imply the result.
\hfill $\blacksquare$ \\

We can transform this result into a result on exponential moments. This is obtained in the exact same way as Lemma 6.2 from Lemma 6.1 in \cite{CJK} and we skip the details. Notice that the modulus of continuity in the hyperbolic case doesn't depend on $\kappa$, unlike Lemmas 6.1 and 6.2 of \cite{CJK}.

\begin{lemma} \label{lem:exp_cont}
Let $u$ be the solution to \eqref{eqn:wave} with $0 < \epsilon_0 := \inf_{z \in \R} \sigma(z) \leq \sup_{z \in \R} \sigma(z) := S_0 < \infty$. Then, for every $t > 0$, there exists a constant $C \in (0,\infty)$ such that,
\begin{equation}
\E\left[\sup_{\substack{x,x'\in I:\\ |x-x'|\leq \delta}} \exp\left(\frac{|u(t,x)-u(t,x')|^2}{C\delta}\right)\right] \leq \frac{2}{\delta},
\end{equation}
uniformly for every $\delta \in (0,1]$ and every interval $I$ of length at most one.
\end{lemma}

We are now ready to prove Theorem \ref{thm:chaos_bdd}. \\

\noindent{\bf Proof of Theorem \ref{thm:chaos_bdd}.}
The lower bound is obtained in a very similar way as in Theorem \ref{thm:chaos_main}. We use Proposition \ref{prop:pr_bdd} instead of Proposition \ref{prop:pr_lb_log}. We can also update Proposition \ref{prop:Picard} in order to consider a moment bound using Proposition \ref{prop:mom_ub_bdd} instead of \ref{prop:mom_ub}. Altogether, \eqref{eqn:main_1} becomes
\begin{equation} \label{eqn:bdd_1}
\P\left\{\max_{1 \leq j \leq N} |u(t,x_j)| < \lambda\right\} \leq \left(1 - c_1 e^{-c_2 (2\lambda)^2}\right)^N + N C^p (\tilde{\mu}_t p)^p e^{-np} \lambda^{-p}.
\end{equation}
Now, we choose the parameters judiciously: we take $\lambda:=p$, $N:=p e^{c_2 p^2}$, $n=\varrho p$, for some constant $\varrho > 2^{1/\gamma} c_2$.
As a consequence, \eqref{eqn:main_1} becomes
\begin{eqnarray}
\P\left\{\max_{1 \leq j \leq N} |u(t,x_j)| < p\right\} & \leq & e^{-c_1 p} + \exp\left(c_2(2p)^2 + \log(p) + p \log(\tilde{\mu}_t) -\varrho p^2\right) \nonumber \\
& \leq & 2 e^{-c_1 p}.
\end{eqnarray}
A similar argument as in Theorem \ref{thm:chaos_main} leads to the lower bound.

As for the upper bound, we follow the same approach as in Theorem 1.2 of \cite{CJK}. Let $R > 0$ be an integer and let $x_j = -R+j$ for $j=1,\ldots,2R$. Then, we can write
\begin{eqnarray}
\P\left\{\sup_{x \in [-R,R]} u(t,x) > 2\alpha(\log R)^{1/2}\right\} & \leq & \P\left\{\max_{1 \leq j \leq 2R} u(t,x_j) > \alpha(\log R)^{1/2}\right\} \label{eqn:bdd_2} \\
&& + \, \P\left\{\max_{1 \leq j \leq 2R} \sup_{x \in (x_j,x_{j+1})} |u(t,x)-u(t,x_j)| > \alpha(\log R)^{1/2}\right\}. \nonumber
\end{eqnarray}
Now, by Proposition \ref{prop:pr_bdd}, we have
\begin{equation} \label{eqn:bdd_3}
\P\left\{\max_{1 \leq j \leq 2R} u(t,x_j) > \alpha(\log R)^{1/2}\right\} \leq 2R \sup_{x\in\R} \P\left\{u(t,x) > \alpha(\log R)^{1/2}\right\} \leq \const \cdot R^{1-c\alpha^2/\max\{\kappa,\overline{v}_0^2 \kappa ^2\}}
\end{equation}
and, by Chebychev's inequality and Lemma \ref{lem:exp_cont} with $\delta=1$,
\begin{equation} \label{eqn:bdd_4}
\P\left\{\max_{1 \leq j \leq 2R} \sup_{x \in (x_j,x_{j+1})} |u(t,x)-u(t,x_j)| > \alpha(\log R)^{1/2}\right\} \leq \const \cdot R^{1 - \alpha^2/C}.
\end{equation}
Hence, replacing \eqref{eqn:bdd_3} and \eqref{eqn:bdd_4} in \eqref{eqn:bdd_2}, we obtain
\begin{equation}
\sum_{R=1}^{\infty} \P\left\{\sup_{x \in [-R,R]} u(t,x) > 2\alpha(\log R)^{1/2}\right\} \leq \sum_{R=1}^{\infty} R^{1-q\alpha^2},
\end{equation}
where
\begin{equation} 
q:=\min\left\{\frac{c}{\max\{\kappa,\overline{v}_0^2 \kappa^2\}},\frac{1}{C}\right\}.
\end{equation}
The sum is convergent provided $\alpha > (2/q)^{1/2}$. Hence, by the Borel-Cantelli Lemma, we have shown that
\begin{equation}
\limsup_{R \rightarrow \infty} \frac{\sup_{x \in [-R,R]} u(t,x)}{(\log R)^{1/2}} \leq \left(\frac{8}{q}\right)^{1/2} \qquad \mbox{a.s.}
\end{equation}
Clearly, $(8/q)^{1/2} \leq \max\{\kappa^{1/2}, \overline{v}_0 \kappa\}/c$ for all $\kappa > \kappa_0$, where $\kappa_0$ is an appropriate constant. A montonicity argument proves the result for non-integer $R$.
\hfill $\blacksquare$ \\

Theorem \ref{thm:chaos_bdd} essentially draws the same conclusion as its parabolic equivalent (Theorem 1.2 in \cite{CJK}), namely that whenever $\sigma$ is bounded away from $0$ and infinity, then the supremum of the solution $u$ behaves as in the case where $\sigma$ is constant; that is, as a Gaussian process: the supremum is of order $(\log R)^{1/2}$ as in the parabolic case. We also notice that the order in $\kappa$ is the same for the upper and lower bound only if $v_0 \equiv 0$. Indeed, whenever the initial derivative does not vanish, it plays a role in the behavior with respect to $\kappa$. If $v_0 \equiv 0$, we get back to a situation similar to the parabolic one. We also notice that the behavior in $\kappa$ is reversed compared to the parabolic case (as in Theorem \ref{thm:chaos_main}): the supremum is increasing in $\kappa$. However, the uniformity in $\kappa$ of the constants doesn't hold for small values of $\kappa$, exactly as in the parabolic case, even though the behavior is reversed. This might look contradictory, but it is not. Indeed, the loss of uniformity in $\kappa$ is due to the modulus of continuity estimate, which is sharp for large values of $\kappa$, but not good when $\kappa$ is small, both in the parabolic and hyperbolic case. It is pretty easy to see this fact in the hyperbolic case by taking a careful look at \eqref{eqn:cont_gamma}. A more careful study of the continuity of $u$ could perhaps lead to more exact results for small values of $\kappa$.


%




\addcontentsline{toc}{section}{References}

\begin{thebibliography}{99}

\bibitem{bertini_cancrini} Bertini L. \& Cancrini N. \emph{The stochastic heat equation: Feynman-Kac formula and intermittence.} Journal of Statistical Physics, Vol. 78, n°5-6 (1995) 1377-1401.
\bibitem{Burkholder} Burkholder D. L. \emph{Martingale transforms.} Annals of Mathematical Statistics. Vol. 37 (1966) 1494-1504.
\bibitem{BDG} Burkholder D. L. \& Davis B. J. \& Gundy R. F. \emph{Integral inequalities for convex functions of operators on martingales.}	In: \emph{Proceedings of the Sixth Berkeley Symposium on Mathematical	Statistics and Probability.} Vol. II, 223-240, University of California Press, Berkeley, California, 1972.
\bibitem{BG} Burkholder D. L. \&  Gundy R. F. \emph{Extrapolation and interpolation of quasi-linear operators on martingales.} Acta Mathematica Vol. 124 (1970), 249-304.
\bibitem{CK} Carlen E. \& Kree P. \emph{$L^p$ estimates for multiple stochastic integrals.} Annals of Probability. Vol. 19, n°1 (1991) 354-368.  
\bibitem{CarmonaNualart} Carmona, R. A. \& Nualart D. \emph{Random nonlinear wave equations: propagation of singularities.} Annals of Probability, Vol. 16, n°2 (1988) 730-751.	
\bibitem{these} Conus D. \emph{The non-linear stochastic wave equation in high dimensions: existence, Hölder-continuity and Itô-Taylor expansion.} Ph.D Thesis n°4265, EPFL Lausanne (2008).
\bibitem{conus_dalang} Conus D. \& Dalang R.C. \emph{The non-linear stochastic wave equation in high dimensions.} Electronic Journal of Probability Vol. 13, 2008. 
\bibitem{CJK} Conus D. \& Joseph M. \& Khoshnevisan D. \emph{On the chaotic character of the stochastic heat equation, before the onset of intermittency.} Annals of Probability (to appear).
\bibitem{CJK2} Conus D. \& Joseph M. \& Khoshnevisan D. \emph{Correlation-length bounds, and estimates for intermittent islands in parabolic SPDEs.} Preprint, 2011.
\bibitem{CJKS} Conus D. \& Joseph M. \& Khoshnevisan D. \& Shiu S.-Y. \emph{On the chaotic character of the stochastic heat equation, II.} Preprint (2011).
\bibitem{conus_khoshnevisan} Conus D. \& Khoshnevisan D. \emph{On the existence and position of the farthest peaks of a family of stochastic heat and wave equations.} Probability Theory and Related Fields (to appear).
\bibitem{dalang} Dalang R.C. \emph{Extending martingale measure stochastic integral with applications to spatially homogeneous spde's.} Electronic Journal of Probability, Vol. 4, 1999.
\bibitem{dalang_frangos} Dalang R.C. \& Frangos N.E. \emph{The stochastic wave equation in two spatial dimensions.} Annals of Probability, Vol. 26, n°1 (1998) 187-212.
\bibitem{dalang_mueller} Dalang R.C. \& Mueller C. \emph{Some non-linear s.p.d.e's that are second order in time.} Electronic Journal of Probability, Vol. 8, 2003.
\bibitem{dalang_mueller_2} Dalang R.C. \& Mueller C. \emph{Intermittency properties in a hyperbolic Anderson problem.} Annales de l'Institut Henri Poincar\'e, Vol. 45, n°4 (2009), 1150-1164.
\bibitem{Davis} Davis B. \emph{On the $L^p$ norms of stochastic integrals and other martingales} Duke Mathematical Journal, Vol.43, n°4 (1976) 697-704.
\bibitem{foondun_khoshnevisan} Foondun M. \& Khoshnevisan D. \emph{Intermittence and nonlinear parabolic stochastic partial differential equations.} Electronic Journal of Probability, Vol. 14, Paper no. 12 (2009) 548--568.
\bibitem{foondun_khoshnevisan_2} Foondun M. \& Khoshnevisan D. \emph{On the global maximum of the solution to a stochastic heat equation with compact-support initial data.} Annales de l'Institut Henri Poincar\'e, Vol. 46, n°4 (2010), 895-907.
\bibitem{hu_nualart} Hu Y. \& Nualart D. \emph{Stochastic heat equation driven by fractional noise and local time.} Probability Theory and Related Fields, Vol. 143, n°1-2 (2009) 285-328.
\bibitem{mueller} Mueller C. \emph{On the support of solutions to the heat equation with noise.} Stochastics and Stochastics Reports, Vol.37, n°4 (1991) 225-245.
\bibitem{revuz-yor} Revuz D. \& Yor M. \emph{Continuous Martingales and Brownian Motion.} Third Edition. Grundlehren der Mathematischen Wissenschaften n°293. Springer-Verlag, Berlin, 2009. 
\bibitem{walsh} Walsh J.B. An Introduction to Stochastic Partial Differential Equations. In : \emph{Ecole d'Eté de Probabilités de St-Flour}, XIV, 1984. Lecture Notes in Mathematics 1180. Springer-Verlag, Berlin, Heidelberg, New-York (1986), 265-439.
\end{thebibliography}
\end{document}